\documentclass[a4paper,twoside,
 11pt]{article}

\usepackage{amssymb,amsfonts,amsmath,xypic}
\newcommand{\sddots}{\!\!\!\!\!{~_{\bf{.}}}\!\cdot~\!\!\!^{^{\bf{.}}}}
\begin{document}
\newtheorem{defin}{~~~~Definition}
\newtheorem{prop}{~~~~Proposition}
\newtheorem{remark}{~~~~Remark}
\newtheorem{cor}{~~~~Corollary}
\newtheorem{theor}{~~~~Theorem}
\newtheorem{lemma}{~~~~Lemma}
\newtheorem{ass}{~~~~Assumption}
\newtheorem{con}{~~~~Conjecture}
\newtheorem{concl}{~~~~Conclusion}
\renewcommand{\theequation}{\thesection.\arabic{equation}}

\title{ On feedback 
 classification of control-affine systems with one and two-dimensional 
inputs} 
\date{}
\author {Andrei Agrachev\thanks{ S.I.S.S.A., Via Beirut 2-4,
34014 Trieste Italy and Steklov Mathematical Institute, 
ul.~Gubkina~8, 117966 Moscow Russia; email: 
agrachev@sissa.it} \and Igor Zelenko \thanks{ S.I.S.S.A., 
Via Beirut 2-4, 34014 Trieste Italy; email: 
zelenko@sissa.it}}\maketitle 

\begin{abstract}
The paper is devoted to the local classification of 
generic control-affine systems on an $n$-dimensional 
manifold with scalar input for any $n\geq 4$ or with two 
inputs for $n=4$ and $n=5$, up to state-feedback 
transformations, preserving the affine structure (in 
$C^\infty$ category for $n=4$ and $C^\omega$ category for 
$n\geq 5$). First using the Poincare series of moduli 
numbers we introduce the intrinsic numbers of functional 
moduli of each prescribed number of variables on which a 
classification problem depends. In order to classify affine 
systems with scalar input we associate with such a system the 
canonical frame by normalizing  some structural functions in a
commutative relation of the vector fields, which define our 
control system. Then, using this canonical frame, we 
introduce the canonical coordinates and find a complete 
system of state-feedback invariants of the system. 
It also gives automatically the micro-local (i.e. local in 
state-input space) classification of the generic non-affine 
$n$-dimensional control system with scalar input for $n\geq 3$ 
(in $C^\infty$ category for $n=3$ and in $C^\omega$ 
category for $n\geq 4$). Further we show how the problem of 
feedback-equivalence of affine systems with two-dimensional input in 
state space of dimensions $4$ and $5$ can be reduced to the 
same problem for affine systems with scalar input. In order to 
make this reduction we distinguish the subsystem of our 
control system, consisting of the directions of all 
extremals in dimension $4$ and all abnormal extremals in 
dimension $5$ of the time optimal problem, defined by the 
original control system. In each classification problem 
under consideration we find the intrinsic numbers of 
functional moduli of each prescribed number of variables 
according to its Poincare series. 
\end{abstract}
\medskip

\vskip .2in


\medskip
\section{Introduction}
\indent \setcounter{equation}{0} 

In the paper for the convenience of the presentation 
all objects are $C^\infty$ without special mentioning, although 
all constructions and some statements remain valid in an obvious way 
also in $C^k$ category 
for an appropriate finite $k$. On the other hand, 
some of the statements are known to be true only in the real analytic 
category, which  will be indicated explicitly.

Let $M$ be an 
$n$-dimensional manifold, $f_0, f_1, \ldots, f_r$ be vector 
fields on $M$, $r<n$. Consider the following control-affine 
system with $r$ inputs on $M$ 
\begin{equation}
\label{affcs} \dot q=f_0(q)+\sum_{k=1}^r u_r\,f_r(q), \quad 
q\in M,\quad u_1\ldots \,u_r\in\mathbb R.\end{equation} 
We say that a system of the type \eqref{affcs} is an 
$(r,n)$ control-affine system. We also assume that at 
given point $q_0$ 
\begin{equation}
\label{indep} \dim {\rm 
span}\bigl(f_0(q_0),f_1(q_0),\ldots, f_r(q_0)\bigr)=r+1. 
\end{equation}
Consider the group $FB_{q_0}$ of state-feedback transformations, 
preserving an affine structure and the point $q_0$, i.e. 
transformations of the type 
\begin{equation}
\label{feedback} 
 \left\{\begin{array} {l} q=\Phi(\tilde q)\\ u=
 {\mathcal B}(\tilde q) \tilde 
u+{\mathcal A}(\tilde q)\\ 
q_0=\Phi(q_0),\end{array}\right.,\quad 
u=\left(\begin{array}{c}u_1\\\vdots\\u_r\end{array}\right), 
\end{equation} where $\Phi$ is a diffeomorphism in a 
neighborhood of $q_0$, ${\mathcal A}(q)\in\mathbb R^r$, 
${\mathcal B}(q)$ is a $r\times r$-matrix, $\det{\mathcal 
B} (q_0)\neq 0$. This group of transformations acts 
naturally on the set of germs of systems of the type 
(\ref{affcs}) and defines the equivalence relation, called 
{\it state-feedback equivalence}. The natural question is 
when two germs of the systems of the type (\ref{affcs}) are 
state-feedback equivalent. 

Let us roughly estimate the "number of parameters" in the 
considered equivalence problem. The set of $r$-dimensional 
affine subspaces in $\mathbb{R}^n$ forms 
$(r+1)(n-r)$-dimensional manifold. Therefore, if the 
coordinates on $M$ are fixed, then the control system of 
the type (\ref{affcs}) can be defined by $(r+1)(n-r)$ 
functions of $n$ variables. The group of the coordinate 
changes on $M$ is parameterized by $n$ functions of $n$ 
variables. So, by a coordinate change one can "normalize", 
in general, only $n$ functions among those $(r+1)(n-r)$ 
functions, defining our control system. Thus we may expect 
that the set of orbits of generic germs of systems 
(\ref{affcs}) at $q_0\in M$ w.r.t. the action of the group 
of transformations of the type (\ref{feedback}) can be 
parameterized by $(r+1)(n-r)-n=r(n-r-1)$ arbitrary germs of 
functions of $n$ variables and a number of germs of 
functions, depending on less than $n$ variables (see the 
next section for the discussion about the number of these 
additional functional invariants). 

According to the last estimate the only case when the 
functional parameters (also called functional moduli) are 
not expected in the parameterization of the set of orbits of 
generic germs of systems (\ref{affcs}) is the case $r=n-1$ 
(i.e., corank $1$ control-affine systems). In this case 
under assumption (\ref{indep}) there is the natural 
one-to-one correspondence between the set of control-affine 
systems, up to feedback transformations, and the set of 
differential $1$-forms on the ambient manifold: to any 
affine system (\ref{affcs}) one can assign a unique 
differential $1$-form $\omega$ such that $\omega(f_i)=0$ 
for $i=1,\ldots n-1$ and $\omega(f_0)=1$. So, the 
state-feedback classification of corank $1$ control-affine 
systems satisfying (\ref{indep}) is equivalent to the 
well-known classification of differential $1$-forms w.r.t. 
the action of the group of diffeomorphisms (see, for 
example, \cite{zhi}, section 3 and Appendix C there ). In 
particular, all germs of (\ref{affcs}) such that the 
underlying vector distribution ${\rm span} 
(f_1,\ldots,f_{n-1})$ is contact for odd $n$ or 
quasi-contact for even $n$ (which is generic assumption) 
are state-feedback equivalent to the control-affine system, 
corresponding to the classical Darboux model (note also 
that in \cite{zhi} normal forms for codimension 1 
singularities are given too). In the context of corank $1$ 
control-affine systems it is not worse also to mention 
 the works \cite{jac}, \cite{resp0}, 
\cite{zhires}, where the case 
$f(q_0)\in {\rm span} \bigl(f_1(q_0),\ldots, 
f_{n-1}(q_0)\bigr)$
was treated. 

The case $r=1$, $n=3$ (the smallest dimensions, when the 
functional parameters appear) was treated in \cite{agr2} 
(section 3, Proposition 3.2 there). In particular, it was 
shown that the set of orbits of generic germs of the 
systems w.r.t. the action of the group of transformations 
of the type (\ref{feedback}) can be parametrized by one 
arbitrary function of $3$ variables, $2$ arbitrary 
functions of two variable , and the discrete invariant from 
the set $\{-1,1\}$. 

\begin{remark}
\label{utochn} {\rm Actually, in Proposition 3.2 of 
\cite{agr2} the two functions of two variables satisfy 
certain conditions on coordinate subspaces of some special 
coordinates, which are canonical up to some reflections, 
but instead of these functions one can take their 
appropriate partial derivatives, which are already 
arbitrary. The functions of the parameterization are 
state-feedback invariants up to some reflections in the 
coordinates.}$\Box$ 
\end{remark}
 
%

 In the present paper we make a classification of generic 
germs of systems of the type (\ref{affcs}), up to 
state-feedback equivalence, in the following cases 
\begin{enumerate}
\item $r=1$, $n=4$;
\item $r=1$, $n\geq 5$ in the real analytic category;
\item
$r=2$, $n=4$; 
\item
$r=2$, $n=5$ in the real analytic category. \end{enumerate} 
\indent

In general, statements of the kind "the classification 
problem depends on the tuple of functional invariants, 
consisting of certain number of functions of each number of 
variables" need to be clarified: although the number of 
functional invariants, depending on the maximal number of 
variables, could be found in the more or less rigorous way 
by counting of dimensions arguments (as was done above), it 
seems that the numbers of functions of each number of 
variables less than the maximal one depend on the way of 
normalization of the objects in the considered 
classification problem. Besides, these numbers could be 
changed rather arbitrary by mixing, combining or separating 
the formal Taylor series of this functional invariants 
without losing any information (at least if we work in the 
category of formal Taylor series or in the real analytic 
category). In \cite{arn} (section 1 there) it was proposed 
to use the so-called {\it Poincare series of the moduli 
numbers of the classification problem} in order to 
determine intrinsically the number of functions of each 
prescribed number of variables, on which some 
classification problem depends. In section 2 below, using 
the Poincare series, 
 we give  a canonical selection of these  numbers. 
The way the canonical parameterization is obtained indicates the 
presence of an interesting algebraic structure on the set of all tuples  
of fundamental invariants parameterizing given classification problem. 
For the moment, this algebraic 
structure remains hidden and needs further research.  

In the case of scalar input  our method of the classification is 
similar to the procedure, used in \cite{agr2} for the case 
$n=3$ and it is described in section 3. It consists 
basically of the following two steps: first for any control 
system, satisfying some genericity assumptions, we 
construct the canonical frame by normalizing some structural 
functions in a commutative relation of the vector 
fields, which define our control-affine systems; then, 
using this canonical frame, we introduce the canonical 
coordinates and find the complete system of state-feedback 
invariants of the system. 

Besides, to any control system 
\begin{equation}
\label{nonaff} \dot y={\mathcal F}(y,v),\quad y\in S,v\in V 
\end{equation} 
on an $m$-dimensional manifold $S$ (the state space) with 
one-dimensional control space $V$ one can assign the 
following control-affine system on the $(m+1)$-dimensional 
state-space $S\times V$. 
\begin{equation}
\label{affcor} \left\{\begin{array}{l} \dot y ={\mathcal F}
(y,v)\\ \dot v=u\end{array}\right., \quad v\in \mathbb{R} 
\end{equation}
(here we look on $u$ as on a new state variable, $v$ is the 
new control, $f_0=({\mathcal F} (y,v),0)^T$ and 
$f_1=(0,1)^T$ in the notations of (\ref{affcs})). It turns 
out that having the local classification of generic 
$m+1$-dimensional control-affine systems with scalar input, 
one gets also the micro-local (i.e. local in state-input 
space) classification of the generic $m$-dimensional 
control system (see, Remark \ref{affrem} at the 
end of section 3). 

Further, in section 4, we show that the problem of 
state-feedback classification of the control-affine systems 
with two-dimensional input in dimensions $4$ and $5$ can be reduced to 
the previous problem for the control-affine systems with 
scalar input in the same dimensions. In order to make this 
reduction we distinguish the subsystem, consisting of the directions of all extremals in 
dimension $4$ and all abnormal extremals in dimension $5$ 
of the time optimal problem, defined by the original 
control system. 

In each classification problem under consideration we find 
the intrinsic numbers of functional moduli of each 
prescribed number of variables according to its Poincare 
series. 

Finally note that the problem, considered here, is 
different from one, considered in the paper \cite{resp}, 
which has the similar title. In the mentioned paper the 
authors study germs of $n$-dimensional control-affine 
system with scalar input at an equilibrium point $q_0$, 
i.e. when $f_0(q_0)\in \{{\mathbb R} f_1(q_0)\}$. Their 
method is a generalization of technique, developed in 
\cite{kang2} and \cite{kang1}, which is similar to 
classical Poincare-Dulac procedure for normalization of 
vector fields near a stationary point. Therefore in the 
method of \cite{resp} it is crucial that $q_0$ is an 
equilibrium point. Here we classify the control-affine 
systems with scalar control near non-equilibrium point, 
which seems very natural in view of the fact that this 
classification, except the case $n=2$, a priori contains 
functional moduli. 
Besides, the feedback invariants, constructed here for 
generic germs could be used also for the problem with 
equilibrium points by passing to the limit.

\section {Poincare series and the intrinsic number of 
functional invariants.} \indent \setcounter{equation}{0} 

We start with some terminology. Let $M$ be a smooth 
manifold. Fix a point $q_0\in M$. Consider a set ${\mathcal 
O}$ of germs at $q_0$ of smooth objects on $M$ (for 
example, Riemannian metrics, vector distributions, 
control-affine systems) such that the group of local
diffeomorphisms ${\rm Diff}_{q_0}$, preserving the point 
$q_0$, acts naturally on it.  This action defines the equivalence relation 
on $\mathcal O$. 

Denote by $J^k({\mathcal O})$ the space of all $k$-jets at 
$q_0$ of objects from the set ${\mathcal O}$. We say that 
the set $\widetilde{\mathcal O}\subset {\mathcal O}$ is 
{\it a generic subset} of ${\mathcal O}$ if there exists an 
integer $k\geq 0$ and a Zariski open set $U$ in 
$J^k({\mathcal O})$ such that $$\widetilde {\mathcal 
O}=\{{\mathfrak b}\in {\mathcal O}: \text{$k$-jet of 
${\mathfrak b}$ belongs to $U$}\}.$$ 

By classification problem on ${\mathcal O}$ we mean the problem 
to find a system of 
fundamental invariants for objects from some generic subset of 
${\mathcal O}$ such that two generic objects  
are equivalent if and only if 
they have the same systems of fundamental invariants. 
Let $\widetilde {\mathcal O}$ be a generic subset of 
${\mathcal O}$, which is invariant w.r.t. the action of the 
group ${\rm Diff}_{q_0}$.

\begin{defin}
\label{funcinv}
A 
mapping $F$ from the set $\widetilde {\mathcal O}$ to the 
set $C^\infty_0(\mathbb R^l, \mathbb R)$ of germs at $0$ of smooth functions in ${\mathbb R}^l$, which 
is invariant w.r.t. the action of the group ${\rm 
Diff}_{q_0}$ on $\widetilde {\mathcal O}$, is called  a 
functional invariant of $l$ variables of a generic subset of objects from 
${\mathcal O}$. 
\end{defin} 
When the object $\mathfrak b\in \widetilde{\mathcal O}$ is fixed, we will mean by the functional invariant also the value of the mapping $F$ at $\mathfrak b$. We will denote this germ of function by the same letter $F$ without special mentioning.

Now let us describe the Poincare series of the moduli 
numbers of the classification problem. The action of the 
group ${\rm Diff}_{q_0}$ induces the action ${\mathcal 
A}_k$ of some finite dimensional Lie group $G_k$ on the 
space $J^k({\mathcal O})$ for any integer $k\geq 0$. So , 
${\mathcal A}_k$ is a mapping from $G_k\times J^k({\mathcal 
O})$ to $J^k({\mathcal O})$. Given any $\mathfrak b \in 
J^k(\mathcal O)$ let ${\mathcal A}_k^{\mathfrak b}$ be the 
mapping from $G_k$ to $J^k({\mathcal O})$ such that 
${\mathcal A}_k^{\mathfrak b}(\cdot)= {\mathcal A}_k 
(\cdot, \mathfrak b)$. Let $e_k$ be the identity of the 
group $G_k$. Set 

\begin{equation}
\label{moduli}
 m(k)=\dim J^k(\mathcal O)-\max_{b\in J^k(\mathcal O)}  {\rm rank}\,d 
 {\mathcal A}_k^{\mathfrak b}(e_k)=\min_{b\in J^k(\mathcal O)}{\rm corank}\, d 
 {\mathcal A}_k^{\mathfrak b}(e_k)
 \end{equation}
 (here ${\rm rank}\, d 
 {\mathcal A}_k^{\mathfrak b}(e_k)$ and ${\rm corank}\,d 
 {\mathcal A}_k^{\mathfrak b}(e_k)$ are  the rank and the corank  of 
 the differential of ${\mathcal A}_k^{\mathfrak b}$ at $e_k$ respectively).
Roughly speaking $m(k)$ is the dimension of the space of 
orbits w.r.t. the last action ${\mathcal A}_k$.
The number $m(k)$ is called the {\it moduli number of the 
$k$-jets}. {\it The Poincare series of the moduli numbers 
of the classification problem} (or shortly {\it the 
Poincare series of the classification problem}) is by 
definition the following function: 
\begin{equation} \label{poincare} 
M(t)=\sum_{k=0}^\infty m(k)t^k 
\end{equation}
\begin{remark} \label{remcal}{\rm Since the integer-valued function ${\mathfrak b}\mapsto 
{\rm rank}\, d 
 {\mathcal A}_k^{\mathfrak b}(e_k)$ takes its maximal value at a Zarisky open set, in (\ref{moduli})
 we can replace $\mathcal O$ by any its generic subset 
$\widetilde{\mathcal O}$.} $\Box$
\end{remark}

The Poincare series could be useful in evaluating of the 
number of the functional invariants of the given 
number of variables, on which the given classification 
problem depends, because of the following well-known fact: 
if one denotes by $j_l(k)$ the dimension of the space $J^k(\mathbb R^l,\mathbb R$) of 
$k$-jets of functions of $l$-variables, then the 
corresponding Poincare series of numbers $j_l(k)$ satisfies 
\begin{equation} 
\label{poinfunc} \sum_{k=0}^\infty 
j_l(k)t^k=\frac{1}{(1-t)^{l+1}} 
\end{equation}
 (here one uses that 
$j_l(k)=\frac{(l+k)!}{l! k!}$). So, if , for example, the 
Poincare series of some classification problem is equal to 
\begin{equation} \label{poincarex} 
M(t)=t^w\sum_{i=1}^n\frac{p_i}{(1-t)^{i+1}}, 
\end{equation}
where all $p_i$ are nonnegative integers, then it is 
natural to conclude that this problem depends on the tuple, 
consisting of $p_i$ functional invariants of $i$ 
variables for each $1\leq i\leq n$, while the parameter $w$ 
(i.e. the order of zero of the Poincare series $M(t)$ at 
$t=0$) is equal to the minimal $k\geq 0$ such that the 
action of the group $G_k$ on the space of $J^k(\mathcal O)$ is 
non-transitive. 

 Till now nothing is known about the form 
of the functions $M(t)$ for general classification 
problems. For example, the following open question is 
stated in \cite{arn}: Is it true that the Poincare series 
of moduli numbers are rational functions in the most 
classification problems?  
In the next sections we will show by direct computations 
that in all classification problems 1-4 listed in 
Introduction it is true and moreover the function $M(t)$ 
has a unique pole at $t=1$. On the other hand, in all 
considered cases (except the cases $r=1$, $n=3$ or $5$) the 
Poincare series has no a representation of the type 
(\ref{poincarex}) with nonnegative $p_i$. Below we give an 
algorithm to extract the number of functional invariants 
from Poincare series also in these cases. 

From now on we will suppose that the Poincare series $M(t)$ 
of the classification problem is a rational function with a 
unique pole at $t=1$. Let $w_0$ be the order of zero of the 
function $M(t)$ at $t=0$. 

\begin{lemma}
\label{replemma} For any integers $w\geq w_0$ and $l$ there 
exist a unique polynomial $R(t)$ with 
\begin{equation}
\label{deg} \deg R(t)< w-w_0
\end{equation}
 and a unique rational function $Q(t)$ 
with the unique pole at $t=1$ such that 
\begin{equation} \label{normrep} 
M(t)=\cfrac{t^{w_0}R(t)}{(1-t)^{l+1}}+t^w Q(t). 
\end{equation} 
\end{lemma}

{\bf Proof.} Let us fix $l\in\mathbb Z$ and prove the 
existence of a representation of the type (\ref{normrep}) 
for any $w\geq w_0$ by induction in $w$. 

If $w=w_0$, then from the condition \eqref{deg} it follows 
that $R(t)\equiv 0$. Then by definition of order of zero 
the function $Q(t)=\frac{1}{t^{w_0}}M(t)$ is rational with 
the unique pole at $t=1$, which implies (\ref{normrep}). 

Now suppose that a representation of the type 
(\ref{normrep}) exists for some $w=\bar w$, $\bar w\geq 
w_0$, and prove its existence for $w=\bar w+1$. For this 
let $Q(t)$ and $R(t)$ be as in the representation 
(\ref{normrep}) for $w=\bar w$. Denote by 
\begin{equation} \label{Q1} 
Q_1(t)=\frac{1}{t}\left(Q(t)-\cfrac{Q(0)}{(1-t)^{l+1}}\right). 
\end{equation} Then by construction $Q_1$ is also the rational 
function with the unique pole at $t=1$. Expressing $Q(t)$ 
from (\ref{Q1}) and substituting it into (\ref{normrep}),
one has 
\begin{equation} \label{normrep1} 
M(t)=\cfrac{t^{w_0}\bigl(R(t)+Q(0)t^{\bar w-w_0}\bigr)
}{(1-t)^{l+1}}+t^{\bar w+1} Q_1(t). 
\end{equation}
Since $\deg \bigl(R(t)+Q(0)t^{\bar w-w_0}\bigr)<\bar 
w-w_0+1$ it implies the existence of a representation 
(\ref{normrep}) also for $w=\bar w+1$. This completes the 
proof by induction of the existence part of the lemma.

Now let us prove the uniqueness part. If there exists 
another representation of $M(t)$ of the type 
(\ref{normrep}) with a polynomial $\bar R(t)$, $\deg \bar 
R(t)<w-w_0$, and a rational function $\bar Q(t)$ instead of 
$R(t)$ and $Q(t)$, then we have the following identity 
$$\bar R(t)-R(t)= t^{w-w_0}(1-t)^{l+1}\Bigl(Q(t)-\bar 
Q(t)\bigr).$$ It implies that the polynomial $\bar 
R(t)-R(t)$ has zero of order not less than $w-w_0$. On the 
other hand, by assumptions $\deg \bigl(R(t)-\bar 
R(t)\bigr)<w-w_0$, which implies that $R(t)\equiv \bar 
R(t)$ and then also $Q(t)\equiv\bar Q(t)$. $\Box$
\medskip 

We will call the representation \eqref{normrep} (with 
$R(t)$ satisfying \eqref{deg}) the {\sl $(w,l)-$ 
representation of the function $M(t)$}. Let $N$ be the 
order of pole of $(1-t)M(t)$ at $t=1$. 

\begin{defin}
\label{nicedef} The $(w,l)$-representation (\ref{normrep}) 
of $M(t)$ with $R(t)$ and $Q(t)$ satisfying 
\begin{equation}
\label{RQ} R(t)=\sum_{i=0}^{w-w_0-1} r_i t^i,\quad 
Q(t)=\sum_{j=l_1}^N \cfrac{q_j}{(1-t)^{j+1}},\quad 
q_{l_1}\neq 0 
\end{equation}
is called nice, if $1\leq l\leq N$, $l_1\geq l$, and all 
coefficients $r_i$, $q_j$ in \eqref{RQ} are nonnegative 
integers. 
\end{defin} 
 
Of course, in general a rational function 
$M(t)=\frac{t^{w_0}Z(t)}{(1-t)^{N+1}}$, where $Z(t)$ is a 
polynomial (even with integer coefficients), may not have 
any nice $(w,l)-$representation. But if the function $M(t)$ 
is the Poincare series of a classification problem, which 
can be parameterized by functional invariants in some 
reasonable way, then $M(t)$ has at least one nice 
representation. 

To be more precise and to explain why the nice 
representation of the Poincare series are interesting let 
us introduce an additional terminology. 
Let $F$ be a functional invariant of $l$ variables
of a generic subset $\widetilde{\mathcal O}$ of objects from
${\mathcal O}$. Denote by $\pi_k$ the natural projection from the set 
$\widetilde{\mathcal O}$ to the space $J^k(\widetilde{\mathcal O})$. 
Let $(x_1,\ldots x_l)$ be the standard coordinates in ${\mathbb R}^l$. 
Let $C^\infty_0(\mathbb R^l,\mathbb R)$ be the set  of germs at $0$ of 
smooth functions in $\mathbb R^l$.
Given a tuple 
$(i_1,\ldots, i_l)$ of nonnegative integers denote by $D^{(i_1,\ldots i_l)}$ 
the mapping from $C^\infty_0(\mathbb R^l,\mathbb R)$ to $\mathbb R$, 
which assigns to any function 
$f:{\mathbb R}^l\mapsto {\mathbb R}$ the value of its partial derivative 
$\frac{\partial^{i_1+\ldots i_l}f}{\partial x_1^{i_1}\ldots \partial x_l^{i_l}}$ at $0$.

\begin{defin}
\label{weightinv} 
We say that the  functional invariant $F$ of $l$ variables of  
a generic subset $\widetilde{\mathcal O}\subset {\mathcal O}$ has 
the weight $w$ if 
for any integer $k\geq w$ the following condition holds: for any tuple 
$(i_1,\ldots, i_l)$ of nonnegative integers with 
$\displaystyle{\sum_{s=1}^l} i_s\leq k-w$
there exists a mapping 
${\mathfrak G}_{(i_1,\ldots i_l)}:J^k(\widetilde{\mathcal  O})\mapsto 
\mathbb R$ such that the following diagram 
\begin{equation*}
\xymatrix{
 \widetilde{\mathcal O} \ar[r]^F \ar[d]_{\pi_k}& C^\infty_0(\mathbb R^l,\mathbb R) \ar[d]^{D^{(i_1,\ldots i_l)}}\\ J^k(\widetilde{\mathcal O})
   \ar[r]_{{\mathfrak G}_{(i_1,\ldots i_l)}}& 
\mathbb R}
 \end{equation*}
is commutative, while for any tuple 
$(i_1,\ldots, i_l)$ of nonnegative integers with $\displaystyle{\sum_{s=1}^l} i_s> k-w$
such mapping ${\mathfrak G}_{i_1,\ldots i_l}$ does not exist.
%
%
\end{defin}
Essentially the weight of the functional invariant $F$ 
is the integer $w\geq 0$ such that for any $k\geq w$ 
the $k$-jets of our objects  depend on all derivatives of $F$ of 
order not greater than $k-w$ but do not depend on derivatives of $F$ 
of order greater 
than $k-w$. 
 

The classification problem with functional moduli will be 
called {\sl regular}, if there exist integers $w_1$ and 
$w_2$, $n_1$, and $n_2$ ( where $0\leq w_1\leq w_2$, $1\leq 
n_1\leq n_2$) and a $(n_2-n_1+1)\times(w_2-w_1+1)$ matrix 
$P$ with nonnegative entries $p_{ij}$ 
such that the set of orbits of generic germs of systems 
(\ref{affcs}) w.r.t. the group ${\rm Diff}_{q_0}$ can be 
parameterized by the set of functional invariants , 
consisting of $p_{ij}$ functional invariants of $i+n_1-1$ 
variables and weight $j+w_1-1$ for any $1\leq i\leq 
n_2-n_1+1$ and $1\leq j\leq w_2-w_1+1$. Without loss of 
generality we suppose also that on the first and the last 
column and row of the matrix $P$ there is at least one 
nonzero entry. 
Directly from the definition of weight, formula (\ref{poinfunc}), and Remark \ref{remcal} it 
follows that the Poincare series of such classification 
problem satisfies 
\begin{equation}
\label{poinreg} M(t)=\cfrac 
{t^{w_1-1}}{(1-t)^{n_1}}\sum_{j=1}^{w_2-w_1+1} 
t^j\left(\sum_{i=1}^{n_2-n_1+1}\frac{p_{ij}}{(1-t)^{i+1}}\right). 
\end{equation}
Note that $w_1=w_0$ and $n_2=N$, where as before $w_0$ is 
the order of zero of $M(t)$ at $t=0$ and $N$ is the order 
of pole of $(1-t) M(t)$ at $t=1$, but all other parameters, 
appearing in (\ref{poinreg}) could not be uniquely 
recovered from the Poincare series $M(t)$. In the 
considered situation we will say that the classification 
problem admits {\it a regular $(w_2, n_1)$-parameterization 
with parameterization matrix $P$}. 

Given some $(w_2, n_1)$-parameterization of the 
classification problem one can easily build new 
$(w_2,n_1)$-parameterization with another parameterization 
matrices.  
Indeed, take some functional invariant $F$ of the weight 
$j_0$, depending on $i_0$ variables, say $x_1,x_2,\ldots, 
x_{i_0}$, where $2\leq i_0\leq n_2-n_1+1$ and $1\leq 
j_0\leq w_2-w_1$. Let $G(x_1,\ldots,x_{i_0})$ be the 
function such that 
\begin{equation}
\label{FG} 
F(x_1,\ldots,x_{i_0})=F(x_1,\ldots,x_{i_0-1},0)+x_{i_0}G(x_1,\ldots, 
x_{i_0}). 
\end{equation}
Then we can obtain the new parameterization of the 
classification problem by replacing the functional 
invariant $F(x_1,\ldots,x_{i_0})$ by two functional 
invariants $F(x_1,\ldots,x_{i_0-1},0)$ and 
$G(x_1,\ldots,x_{i_0})$. Obviously, the first invariant has 
weight $j_0$ and depends on $i_0-1$ variables, while the 
second one has weight $j_0+1$ and depends on $i_0$ 
invariants. The matrix of the new parameterization is obtained 
from the original one by decreasing the $(i_0,j_0)$-entry 
by 1 and increasing both $(i_0-1,j_0)$-entry and 
$(i_0,j_0+1)$ entry by 1. Such transformation on the set of 
$(N-n_1+1)\times(w_2-w_0+1)$ matrices will be called {\it 
an elementary transformation}. Conversely, given two 
functional invariants $G_1$ and $G_2$ such that $G_1$ 
depends on $i_0-1$ variables say $x_1,\ldots,x_{i_0-1}$ and 
has the weight $j_0$, while $G_2$ depends on $i_0$ 
variables say $x_1,\ldots,x_{i_0}$ and has the weight 
$j_0+1$ (here again 
 $2\leq i_0\leq n_2-n_1+1$ and $1\leq 
j_0\leq w_2-w_1$) one can build the new parameterization by 
replacing the invariants $G_1$ and $G_2$ by one invariant
$G_1+x_{i_0} G_2$, which depends on $i_0$ variables and has 
the weight $j_0$. Of course, in this case the matrix of the 
new parameterization is obtain from the original one by the 
transformation, which is inverse to the elementary one. 
%
%

Now for convenience denote 
\begin{equation}
\label{dim} K_1= N-n_1+1,\quad K_2=w_2-w_0+1.\end{equation} 
Note that among all matrices, which can be obtained from 
the given $K_1\times K_2$ matrix $P$ with integer entries 
by a composition of a finite number of elementary 
transformations and their inverses, there exists a unique 
matrix, denoted by ${\rm Norm}(P)$, such that all its 
entries, except those lying on the first row and the last 
column, are equal to zero. To prove the existence of ${\rm 
Norm}(P)$ one can vanish the entries of the matrix $P$ by a 
composition of elementary transformations and their 
inverses step by step, starting from the entry in the 
left-lower corner, going along the first column from the 
bottom to the top till the entry on the second row , then 
passing to the bottom of the second column, going along it 
from the bottom to the top till the entry on the second row 
and so on till the column before the last one. The 
uniqueness follows from the fact that if we put the entries 
of the matrix ${\rm Norm}(P)$ instead of the entries of $P$ 
in the representation \eqref{poinreg}, then we obtain the 
$(w_2, n_1)$-representation of the Poincare function 
$M(t)$. This representation is unique according to Lemma 
\ref{replemma} and the matrix ${\rm Norm}(P)$ is obviously 
uniquely recovered from it. Also it is not difficult to 
express all nontrivial entries of ${\rm Norm}(P)$ by the 
entries of $P$: 
\begin{subequations} \label{norm}
\begin{align}
~& \left({\rm 
Norm}(P)\right)_{1j}=p_{1j}+\sum_{l=0}^{j-1}\sum_{k=1}^{K_1-1} 
\binom{k+l-1}{l}p_{k+1,j-l},&~&\quad 1\leq j\leq 
K_2-1,\label{norm:a}\\ 
 ~& \left({\rm 
Norm}(P)\right)_{i,K_2}=p_{i,K_2}+\sum_{l=0}^{K_1-i+1}\sum_{k=1}^{K_2-1} 
\binom{K_2-k+l-1}{l}p_{i+l,k},&~&\quad 2\leq i\leq 
K_1,\label{norm:b}\\ ~& \left({\rm Norm}(P)\right)_{1,K_2} 
=p_{1,K_2}&~&\label{norm:c} 
\end{align}
\end{subequations} 
The last relation can be proved, for example, using the 
procedure of passing from $P$ to ${\rm Norm}(P)$, described 
above, and the following well-known combinatorial identity: 
$$\sum_{i=1}^n\binom{i+k-1}{k}=\binom{n+k}{k+1}.$$ 

If the matrix $P$ has only nonnegative integer entries, 
then the matrix ${\rm Norm}(P)$ is obtained from $P$ by a 
finite composition of elementary transformations (without 
using their inverses) and has also only nonnegative integer 
entries (which follows also from relations \eqref{norm}). 
Moreover, if we put the entries of the matrix ${\rm 
Norm(P)}$ instead of the entries of $P$ in the 
representation \eqref{poinreg}, then we obtain the nice 
$(w_2, n_1)$-representation of the Poincare function $M(t)$ 
of our classification problem. We also say that this 
nice representation {\it corresponds to the matrix ${\rm 
Norm}(P)$}. We can summarize all above in the following 

\begin{prop}
\label{summprop} If the classification problem admits a 
regular $(w_2, n_1)$-parameterization with parameterization 
matrix $P$, then it admits a regular $(w_2, 
n_1)$-parameterization with parameterization matrix ${\rm 
Norm}(P)$ and its Poincare series has the nice 
$(w_2,n_1)$-representation, which corresponds to the matrix 
${\rm Norm}(P)$. 
\end{prop} 

The last proposition indicates that the nice 
representations of the Poincare series (if they exist) may 
be used in the definition of the intrinsic number of 
functional invariants of each number of variables and 
weight, on which the given classification problem depends. 
Suppose that the Poincare series $M(t)$ has the nice 
representation for some $(w,l)$. So, the set 
\begin{equation}
\label{niceset} {\rm NS} 
\bigl(M(t)\bigr)\stackrel{def}{=}\{(w,l): {\rm the 
}\,(w,l){\rm -representation}\,\,{\rm of}\,\, M(t)\,\,{\rm 
is}\,\,{\rm nice}\} 
\end{equation}
is not empty. The natural question is what pair to choose 
from ${\rm NS} \bigl(M(t)\bigr)$? To answer this question 
we propose to introduce the order $\prec$ on the set of 
ordered pairs $(w,l)$ in the following way: 
$(w,l)\prec(\bar w,\bar l)$ if and only if $w<\bar w$ or 
$w=\bar w$, but $l>\bar l$. By Definition \ref{nicedef} 

\begin{equation}
\label{Nicein1} {\rm NS} \bigl(M(t)\bigr)\subset\{(w,l): 
w\geq w_0, l\leq N\}, \end{equation} which implies 
immediately that the set ${\rm NS} \bigl(M(t)\bigr)$ 
contains the minimal element w.r.t. the introduced order 
$\prec$. This minimal element will be called {\it the 
characteristic pair of the classification problem}. Denote 
it by  
$(\bar w,\bar l)$. 
Let ${\mathcal C}$ be the $(N-\bar l +1)\times(\bar 
w-w_0+1)$ matrix such that the $(\bar w,\bar 
l)$-representation of $M(t)$ corresponds to the matrix 
${\mathcal C}$.
\begin{defin}
\label{intrinsdef} The $(ij)$-entry of the matrix 
${\mathcal C}$ is called the intrinsic number of the 
functional invariants of $i+\bar l-1$ variables and the 
weight $j+w_0-1$ of the considered classification problem. 
The matrix ${\mathcal C}$ is called the characteristic 
matrix of the classification problem. Any regular $(\bar 
w,\bar l)$-parameterization of the problem (if it exists) 
with the parameterization matrix ${\mathcal C}$ is called 
the characteristic regular parameterization. 
\end{defin} 
In general, it is better to have  a parameterization, consisting of invariants, which  
have  minimal possible weight and depend on maximal possible number of 
variables.
Our definition of characteristic parameterization is in accordance with 
this goal.  
Actually the maximal weight of invariants, appearing in  
characteristic regular parameterization is not greater than the maximal 
weight of invariants, appearing in any other regular parameterization. 
Besides the minimal number of variables in  invariants of regular parameterization is not less than the minimal number of variables in invariants of any other regular parameterization, having the same maximal weight of invariants as a characteristic one. 

One can improve the formula \eqref{Nicein1} for 
the localization of the set ${\rm NS}\bigl(M(t)\bigr)$. 
Indeed let $d$ be the degree of the rational function 
$M(t)$ at infinity. Namely, if 
$M(t)=\frac{Q_1(t)}{Q_2(t)}$, where $Q_1(t)$ and $Q_2(t)$ 
are polynomials, then $d= \deg Q_1(t)-\deg Q_2(t)$. Then 
\begin{equation}
\label{Nicein2} {\rm NS} \bigl(M(t)\bigr)\subset\{(w,l): 
w\geq w_0, 1\leq l\leq \min\,(w-d-1,N)\}. \end{equation} To 
prove (\ref{Nicein2}) we actually have to prove that if the 
pair $(w,l)\in {\rm NS} \bigl(M(t)\bigr)$ then $l\leq 
w-d-1$. Indeed, from (\ref{normrep}) and (\ref{RQ}) it 
follows that $d=\max\, (w-l_1-1,w_0+\deg R-l_1-1)$. But 
first, since $l_1>l$, we have $w-l_1-1\leq w-l-1$ and 
secondly, since $\deg R<w-w_0$, we have $w_0+\deg 
R-l_1-1<w-l-1$. Therefore $d\leq w-l-1$, Q.E.D. 

From (\ref{Nicein2}) it follows also that
\begin{equation}
\label{Nicein3} {\rm NS} \bigl(M(t)\bigr)\subset\{(w,l): 
w\geq \max\, (w_0, d+2) \}. 
\end{equation} 
 The relations (\ref{Nicein2}) and (\ref{Nicein3}) may be useful in 
searching for the characteristic pair of the classification 
problem. Another useful property of the set ${\rm 
NS}\bigl(M(t)\bigr)$ can be formulated as follows: 

\begin{lemma}
\label{netspuska} Assume that the function $M(t)$ has the 
nice $(w,l)$-representation (\ref{normrep}), the functions 
$R(t)$, $Q(t)$, and the number $l_1$ are as in (\ref{RQ}),  
and $l_1=l$ (or, equivalent, $q_l>0$), then 
$(w-1,l-1)\not\in {\rm NS}\bigl(M(t)\bigr)$. 
\end{lemma} 

{\bf Proof.} Let $S(t)$ be the polynomial such that 
$M(t)=\frac{S(t)}{(1-t)^{N+1}}$. Then, using the assumption 
$l=l_1$, it is easy to get 
\begin{equation} \label{degS} \deg S(t)=w+N-l.
\end{equation}
Moreover, directly from (\ref{normrep}) and (\ref{RQ}) one 
can obtain that \begin{equation} \label{leadcoef} 
\cfrac{d^{w+N-l}S}{d\,t^{w+N-l}}=(-1)^{N-l}(w+N-l)!\,q_l. 
\end{equation} 
On the other hand, if the $(w-1,l-1)$-representation of 
$M(t)$ has the form 
\begin{equation}
\label{nonnice} M(t)=\cfrac{t^{w_0}\sum_{i=0}^{w-w_0-2} 
\bar r_i t^i }{(1-t)^{l}}+t^{ w-1} \sum_{j=l_2}^N 
\cfrac{\bar q_j}{(1-t)^{j+1}},\quad \bar q_{l_2}\neq 0, 
\end{equation}
 then $\deg S(t)=\max(w-1+N-l_2, w+N-l-1)$. Comparing this 
with (\ref{degS}) one gets easily that $l_2=l-1$. But then 
by analogy with (\ref{leadcoef}) (applied for 
$(w-1,l-1)$-representation instead of 
$(w,l)$-representation) one has \begin{equation} 
 \label{leadcoef1} 
\cfrac{d^{w+N-l}S}{d\,t^{w+N-l}}=(-1)^{N-l+1}(w+N-l)!\,\bar 
q_{l-1}. 
\end{equation} 
Comparing (\ref{leadcoef}) and (\ref{leadcoef1}), we obtain 
that $\bar q_{l-1}=-q_l$. Hence $\bar q_{l-1}<0$ and the 
$(w-1,l-1)$-representation (\ref{nonnice}) is not nice. 
$\Box$ \medskip 

As a direct consequence of Proposition \ref{summprop}, the 
previous lemma, and the relation (\ref{norm:c}) one has the 
following 
\begin{cor}
\label{usecor} If the classification problem with the 
Poincare series $M(t)$ admits the regular 
$(w,l)$-parameterization with the parameterization matrix $P$ 
such that the entry in the right-upper corner of $P$ is 
positive (i.e., in the previous notations $p_{1, 
w-w_0+1}>0$), then $(w-1,l-1)\not\in {\rm 
NS}\bigl(M(t)\bigr)$. 
\end{cor} 

In the sequel we will show that all classification problems 
1-4 listed in Introduction are regular. For each of this 
problems we will describe explicitly some its regular 
$(w,l)$-parameterization such that $(w,l)$ is the 
characteristic pair of the problem and find its 
characteristic matrix, which also gives the way to obtain 
the characteristic parameterization. 
 
\begin{remark}
\label{13rem} {\rm Sometimes (as in the case of 
(1,3) control-affine systems) in the parameterization of a 
classification problem appear invariants from the finite 
set in addition to the functional parameters. Also these 
functional parameters may be not invariants by themselves, 
but up to some finite group of transformation of 
coordinates. We call such classification problems {\it 
quasiregular}. Obviously, the additional invariants from 
the finite set and the additional freedom in choosing the 
coordinates, up to the finite group of transformation, do 
not affect on the Poincare series. Therefore the identity 
(\ref{poinreg}) for the Poincare series is true also for 
quasiregular problems, where, as before, $p_{ij}$ is the 
number of functional invariants of $i+n_1-1$ variables and 
weight $j+w_1-1$ for any $1\leq i\leq n_2-n_1+1$ and $1\leq 
j\leq w_2-w_1+1$, and we can extend to the quasiregular 
case the same terminology and the constructions, as in the 
regular one. In particular, according to Proposition 3.2 of 
\cite{agr2} (see also Remark \ref{utochn} and the paragraph 
before it in Introduction), the state-feedback 
classification problem for $(1,3)$ control-affine system 
admits quasiregular $(2,2)$-parameterization with the 
following $2\times 2$ parameterization matrix $P= 
\begin{pmatrix} 2&0\\0&1\end{pmatrix}.$
Besides, $w_0=1$ and $N=3$. By the inverse to the 
elementary transformation one can transform $P$ to the 
matrix 
$\widetilde P=\begin{pmatrix} 1&0\\1&0\end{pmatrix}$. The 
Poincare series $M(t)$ of the problem satisfies 
$$M(t)=t\left(\cfrac{1}{(1-t)^3}+\cfrac{1}{(1-t)^4}\right).$$ 
It is not difficult to see that by erasing the last column 
of the matrix $\tilde P$ one obtains the characteristic 
matrix ${\mathcal C}=(1,1)^T$ of the considered 
classification problem and the characteristic pair is equal 
to $(1,2)$. So, the characteristic quasiregular 
parameterization consists of one function of 3 variables and 
the weight $1$, one function of 2 variables and the weight 
$1$, and the discrete invariant from the set $\{-1,1\}$. 
This parameterization is obtained from the original one by a 
rearrangement of the invariants, which corresponds to the 
inverse to elementary transformation, transforming the 
matrix $P$ to the matrix $\widetilde P$ (such 
rearrangements were described in the paragraph after the 
formula (\ref{FG})).$\Box$ } 
\end{remark} 
 
%


\section
{Classification of $(1,n)$ control-affine systems for 
$n\geq 4$} \indent \setcounter{equation}{0} 

For $r=1$ the system (\ref{affcs}) has the following form 
\begin{equation}
\label{affcs1} \dot q=f_0(q)+u\,f_1(q), \quad q\in M,\quad 
u\in {\mathbb R}.\end{equation} Our genericity assumptions
are \begin{eqnarray} \label{gen1} &~&\dim {\rm span} 
(f_0,f_1,[f_1, f_0],\ldots ({\rm ad} f_1)^{n-2}f_0)=n,\\ 
\label{gen2}&~&\dim {\rm span} (f_1,\bigl[f_0,[f_0, 
f_1]\bigr],[f_1, f_0],\ldots ({\rm ad} f_1)^{n-2}f_0)=n, 
\end{eqnarray} and for $n\geq 5$ also 
\begin{equation}
\label{gen3} \dim {\rm span} (f_0,f_1,\bigl[f_0,[f_0, 
f_1]\bigr],[f_1,f_0],\ldots ({\rm ad} f_1)^{n-3}f_0)=n. 
\end{equation} 

The group of feedback transformations 
\begin{equation}
\label{feedback1}
 u= 
 \beta(q) \tilde 
u+\alpha(q), \beta(q)\neq 0 \end{equation} acts naturally 
on the set of pairs of vector fields $(f_0,f_1)$. 
The orbit w.r.t. 
this action is 
\begin{equation}
\label{orb} {\mathcal O}_{(f_0,f_1)}=\{(f_0+\alpha f_1, 
\beta f_1):\alpha, \beta:M\mapsto \mathbb R\,\,{\rm 
are}\,\, {\rm functions}, \beta\neq 0\}. 
\end{equation} The first observation is given by the following 
\begin{prop}
\label{step1}
If the pair $(f_0,f_1)$ satisfies conditions (\ref{gen1}) 
and (\ref{gen2}), then there exists a unique pair $(F
_0,F_1)\in {\mathcal O}_{(f_0,f_1)}$ such that 
\begin{equation} 
\label{N1} \bigl[F_0,[F_0,F_1]\bigr]=F_0+I_1 F_1+ I_2 
[F_1,F_0]+\sum_{k=3}^{n-2} I_k ({\rm ad }F_1)^k F_0 
\end{equation} 
\end{prop}

{\bf Proof.} By assumptions (\ref{gen1}) the vector fields 
$f_0,f_1,[f_1, f_0],\ldots ({\rm ad} f_1)^{n-2}f_0$ 
constitute the frame on $M$. Therefore there are functions 
$E$, $N$, $J_1,\ldots, J_{n-2}$ such that 
\begin{equation}
\label{prim} \bigl[f_0,[f_0,f_1]\bigr]=E f_0+J_1 f_1+
 J_2 [f_1,f_0]+Z \bigl[f_1,[f_1,f_0]\bigr]+\sum_{k=3}^{n-2} J_k 
({\rm ad} f_1)^k f_0 
\end{equation}
Take some pair $(\tilde f_0,\tilde f_1)\in {\mathcal 
O}_{(f_0,f_1)}$,

\begin{equation}
\label{ab} \tilde f_0=f_0+\alpha f_1,\quad \tilde f_1=\beta 
f_1. 
\end{equation}
Suppose that
 \begin{equation}
\label{primtil} \bigl[\tilde f_0,[\tilde f_0,\tilde 
f_1]\bigr]=\widetilde E \tilde f_0+\widetilde J_1 \tilde 
f_1+ 
 \widetilde J_2 [\tilde f_1,\tilde f_0]+\widetilde Z \bigl[\tilde f_1,
 [\tilde f_1,\tilde f_0]\bigr]+\sum_{k=3}^{n-2} J_k 
({\rm ad} \tilde f_1)^k \tilde f_0 
\end{equation}
First note that 
\begin{equation}
\label{transE} \widetilde E=\beta E. \end{equation} It 
follows immediately from (\ref{ab}) and the following 
relations \begin{eqnarray} &~&\label{rel1} \bigl[\tilde 
f_0,[\tilde f_0,\tilde 
f_1]]\bigr]\equiv\beta\bigl[f_0,[f_0, 
f_1]\bigr]-\alpha\beta\bigl[f_1,[f_1,f_0]\bigr]\,\,{\rm 
mod}\,{\rm span}(f_1,[f_0,f_1]),\\&~&\label{rel2} ({\rm ad} 
 f_1)^k  f_0 \in{\rm span}\bigl(\tilde 
f_1,\ldots,({\rm ad} \tilde f_1)^k \tilde f_0 \bigr),\,\,\, 
k\in \mathbb N. 
\end{eqnarray}
From assumption (\ref{gen2}) it follows that $E\neq 0$. 
Therefore, taking $\beta=\frac{1}{E}$, we make 
\begin{equation}
\label{norm1} \widetilde E= 1. 
\end{equation}
Let us denote by $\overline {\mathcal O}_{(f_0,f_1)}$ the 
set of all pairs $(\tilde f_0,\tilde f_1)$, satisfying 
(\ref{norm1}). 
We can assume from the beginning that the original pair 
$(f_0,f_1)$ belongs to $\overline {\mathcal 
O}_{(f_0,f_1)}$, i.e. $E=1$ (we make this assumption just 
in order to avoid extra notations). If $(\tilde f_0,\tilde 
f_1)\in \overline {\mathcal O}_{(f_0,f_1)}$, then also 
$\tilde E=1$. Hence $\beta=1$ or equivalently
$f_1=\tilde f_1$. 
 In other 
words, condition (\ref{norm1}) normalizes the vector field 
$f_1$ or the direction defining the straight line in the 
set of admissible velocities of the system (\ref{affcs1}) 
at any point.  

Further, from (\ref{ab}), taking into account that 
$\beta=1$, it follows easily that $$ ({\rm ad} 
 f_1)^k  f_0 \equiv ({\rm ad} \tilde f_1)^k \tilde f_0\,\,{\rm 
mod}\,{\rm span}(f_1),\,\,\, k\in \mathbb N.$$ This and 
relation (\ref{rel1}) imply that \begin{equation} 
\label{transZ} \widetilde Z=Z-\alpha. \end{equation} 
Setting $\alpha=Z$, we make $\tilde Z=0$, which normalizes 
the drift $\tilde f_0$. So, we have proved that there is a 
unique $(\tilde f_0,\tilde f_1)\in {\mathcal 
O}_{(f_0,f_1)}$ such that $\tilde E=1$ and $\tilde Z=0$, 
which completes the proof of the proposition. $\Box$ 

\begin{remark}
{\rm The mappings  $I_1,\ldots, I_{n-2}$ from $M$ to $\mathbb R$, 
defined by identity (\ref{N1}), 
are state-feedback invariants of the control system (\ref{affcs1}).} $\Box$
\end{remark}

 The vector field $F_0$  and the pair of vector fields $(F_0,F_1)$ from  
Proposition \ref{step1} are called {\it the canonical 
drift} and {\it the canonical pair} of the system 
(\ref{affcs1}) respectively. 

\begin{remark}
\label{F0mean}
{\rm 
Actually, in the case $n=4$, the vector $F_0(q)$ is the velocity of the unique 
abnormal extremal starting at $q$ of the time optimal problem defined by 
system (\ref{affcs1}).
} $\Box$ 
\end{remark}

Now fix some point $q_0\in M$. Denote by $e^{tf}$ the flow, 
generated by the vector field $f$ and $q\circ e^{tf}$ the 
image of the point $q$ w.r.t. this flow. Let 
$\Phi_n:\mathbb R^n\mapsto M$ be the following mapping 
\begin{equation}
\label{F4}\Phi_4(x_1,x_2,x_3,x_4)=q_0\circ e^{x_4 
\bigl[F_1,[F_1,F_0] \bigr]}\circ e^{x_3 [F_1,F_0]}\circ 
e^{x_2 F_1}\circ e^{x_1 F_0},
\end{equation}
\begin{equation}
\label{Fn}
\begin{array}{lll}
\Phi_n(x_1,\ldots,x_n)=&q_0\circ e^{x_n({\rm ad} 
F_1)^{n-3}F_0}\circ\ldots \circ 
e^{x_5\bigl[F_1,[F_1,F_0]\bigr]}\circ\\ ~&\circ 
e^{x_4\bigl[F_0,[F_0,F_1]\bigr]}\circ e^{x_3[F_1,F_0]}\circ 
e^{x_2 F_0}\circ e^{x_1 F_1}, \quad n\geq 5 
\end{array}
\end{equation} 
From assumption (\ref{gen1}) in the case $n=4$ or 
assumption (\ref{gen3}) in the case $n\geq 5$ it follows 
that $\Phi_n'(0)$ is bijective. Hence $\Phi_n^{-1}$ defines 
{\it the canonical coordinates in a neighborhood of $q_0$} 
(or shortly {\it the canonical coordinates at $q_0$}). 
Denote 
\begin{equation}
\label{cI} {\mathcal I}_k=I_k\circ\Phi_n,\quad k=1,\ldots 
n-2. \end{equation} 
Assigning to any generic germ at $q_0$ of  
control-affine  systems (\ref{affcs1}) the function $ {\mathcal I}_k$ we obtain the functional invariant of $n$ variables of this set of objects in the sense of Definition \ref{funcinv} for any $1\leq k\leq n-2$.

Now let us consider the cases $n=4$ and 
$n\geq 5$ separately: 

{\bf a) The case $n=4$.} By (\ref{F4}) and (\ref{N1}), in 
the canonical coordinates the vector fields $F_0$ and $F_1$ 
have the following form: \begin{equation} \label{nf4} 
F_0=\frac{\partial}{\partial x_1},\quad F_1=\sum_{k=1}^4 
a_k \frac{\partial}{\partial x_k}, \end{equation} where the 
components of $F_1$ satisfy the following second order 
linear ordinary differential equations 
w.r.t. the variable $x_1$ 
\begin{equation}
\label{condnf4} 
\frac{\partial^2 a_k}{\partial x_1^2}+{\mathcal 
I}_2\frac{\partial a_k}{\partial x_1}-{\mathcal I}_1 
a_k-\delta_{1,k}=0\quad k=1,2,3,4, 
\end{equation}
with the following restrictions on the initial conditions 
for any $k=1,2,3,4$ 
\begin{subequations}\label{bound4}
\begin{gather}
a_k(0,x_2,x_3,x_4)\equiv \delta_{2k} \label{bound4:a}, 
\\\cfrac{\partial 
a_k}{\partial x_1}(0,0,x_3,x_4)
\equiv -\delta_{3k}\label{bound4:b},\\ \cfrac{\partial^2 
a_k}{\partial x_1\partial x_2}(0,0,0,x_4) 
\equiv -\delta_{4k} \label{bound4:c}, 
\end{gather}
 \end{subequations} 
where $\delta_{ij}$ is the Kronecker symbol. Let for any 
$k=1,2,3,4$ 
\begin{subequations}
\label{prime4} \begin{gather} 
\beta_k(x_2,x_3,x_4)\stackrel{def}{=}\cfrac{\partial^3 
a_k}{\partial x_1 \partial 
x_2^2}(0,x_2,x_3,x_4),\label{prime4:a}\\ 
\psi_k(x_3,x_4)\stackrel{def}{=} \cfrac{\partial^3 
a_k}{\partial x_1 
\partial x_2\partial x_3}(0,0,x_3,x_4)\label{prime4:b}. 
\end{gather}
\end{subequations} 
So, with any germ at $q_0$ of a four-dimensional affine 
system (\ref{affcs1}), satisfying genericity assumptions 
(\ref{gen1}) and (\ref{gen2}), one can associate the 
ordered tuple 
\begin{equation}
\label{prince4} ({\mathcal I}_1, {\mathcal 
I}_2,\beta_1,\beta_2,\beta_3,\beta_4,\psi_1,\psi_2,\psi_3,\psi_4), 
\end{equation}
of state-feedback functional invariants, consisting of two 
germs ${\mathcal I}_1$ and ${\mathcal I}_2$ of functions of 
four variables at $0$, four germs 
$\beta_1,\beta_2,\beta_3,\beta_4$ of functions of three 
variables at $0$, and four germs 
$\psi_1,\psi_2,\psi_3,\psi_4$ of functions of three 
variables at $0$. We call it {\it the tuple of the primary 
invariants of the $(1,4)$ control-affine system} 
(\ref{affcs1}) {\it at the point $q_0$}. Note that by 
(\ref{prime4}) the functional invariants $\beta_k$ and 
$\psi_k$ have the weight $3$ for any $1\leq k\leq 4$, while 
by \eqref{condnf4} and \eqref{bound4} the functional 
invariants ${\mathcal I}_1$ and ${\mathcal I}_2$ have the 
weight $2$. 

Further, fixing $\beta_k$ and $\psi_k$ and using 
\eqref{bound4:b} and \eqref{bound4:c}, one can find 
$\frac{\partial a_k}{\partial x_1}(0,x_2,x_3,x_4)$ for any 
$1\leq k\leq 4$ by the appropriate integrations (see 
(\ref{BK}) below). If in turn one fixes also ${\mathcal 
I}_1$ and ${\mathcal I}_2$, then from the knowledge of 
$\frac{\partial a_k}{\partial x_1}(0,x_2,x_3,x_4)$, 
condition \eqref{bound4:a}, and differential equation 
(\ref{condnf4}) we can recover the functions 
$a_k(x_1,x_2,x_3,x_4)$ and therefore our control-affine 
system itself, just using the standard existence and 
uniqueness results from the theory of ordinary differential 
equations. We summarize all above in the following: 

\begin{theor}
\label{class4}
Given two arbitrary germs ${\mathcal I}_1$ and ${\mathcal 
I}_2$ of functions of four variables at $0$, four arbitrary 
germs $\beta_1,\beta_2,\beta_3,\beta_4$ of functions of 
three variables at $0$, and four arbitrary germs 
$\psi_1,\psi_2,\psi_3,\psi_4$ of functions of two variables 
at $0$ there exists a unique, up to state-feedback 
transformation of the type (\ref{feedback}), 
four-dimensional control-affine system with scalar input, 
satisfying genericity assumptions (\ref{gen1}) and 
(\ref{gen2}), such that the tuple $({\mathcal I}_1, 
{\mathcal I}_2,\beta_1,\beta_2,\beta_3,\beta_4,\psi_1, 
\psi_2,\psi_3,\psi_4)$ is its tuple of the primary 
invariants at the given point $q_0$. In other words, the 
tuples of the primary invariants give the regular 
$(3,2)$-parameterization of the considered classification 
problem with the following $3\times 2$ parameterization 
matrix $P$: 
\begin{equation}
\label{parmat4} P=\begin{pmatrix} 0&4\\ 
0&4\\2&0\end{pmatrix}. \end{equation} The Poincare series 
$M(t)$ of the considered classification problem satisfies 
\begin{equation}
\label{poin4} M(t)=\cfrac{2t^2}{(1-t)^5}+t^3 
\left(\cfrac{4}{(1-t)^4}+\cfrac{4}{(1-t)^3}\right). 
\end{equation}
\end{theor} 

It turns out that $(3,2)$ is the characteristic pair of the 
considered classification problem. Indeed, let as before 
$w_0$ be the order of zero of $M(t)$ at $t=0$, $N$ is the 
order of pole of $(1-t)M(t)$ at $t=1$, and $d$ is the 
degree of $M(t)$ (at infinity). Then from (\ref{poin4}) it 
follows that $w_0=2$, $N=4$, and $d=0$. Hence from 
(\ref{Nicein1}) (or \eqref{Nicein3}) it follows that 
\begin{equation}\label{NS41}
{\rm NS}\bigl(M(t)\bigr)\subset\{(w,l):w\geq 2\} .
\end{equation}
Further, using \eqref{Nicein2}, one has
\begin{equation}\label{NS42}
{\rm NS}\bigl(M(t)\bigr)\cap\{(w,l):w= 
2\}\subset\{(2,1)\},\quad{\rm 
NS}\bigl(M(t)\bigr)\cap\{(w,l):w= 3\}\subset\{(3,1), 
(3,2)\} . 
\end{equation} 
But by the previous theorem our classification problem 
admits $(3,2)$-parameterization such that its 
parameterization metric has a positive entry in the 
right-upper corner. Therefore from Corollary \ref{usecor} 
it follows that $(2,1)\notin {\rm NS}\bigl(M(t)\bigr)$. 
Since $(3,2)\prec(3,1)$ we can conclude from (\ref{NS41}) 
and (\ref{NS42}) that $(3,2)$ is the minimal element of 
${\rm NS}\bigl(M(t)\bigr)$. In other words, $(3,2)$ is the 
characteristic pair of our classification problem. Also, it 
implies that the characteristic matrix ${\mathcal C}$ of 
the problem is equal to ${\rm Norm}(P)$, which can be found 
easily by the series of elementary transformations. Namely, 
\begin{equation}
\label{C4} {\mathcal C}={\rm Norm}(P)=\begin{pmatrix} 2&4\\ 
0&6\\0&2\end{pmatrix}. 
\end{equation}
\begin{concl}
\label{concl14} The characteristic parameterization of 
$(1,4)$ control-affine systems, up to state-feedback 
transformations, consists of $2$ functional invariants of 
$4$ variables and the weight $3$, $6$ functional invariants 
of $3$ variables and the weight $3$, $2$ functional 
invariants of $2$ variables and the weight $2$, and $4$ 
functional invariants of $2$ variables and the weight $3$.
\end{concl} 

In order to obtain a characteristic parameterization from 
the parameterization by the tuple of the primary invariants 
one can implement some series of rearrangement of the 
primary invariants according to the series of elementary 
transformations from the matrix $P$ to ${\rm Norm}(P)$, as 
was described in the previous section (see, for example, 
formula (\ref{FG}) and the paragraph after it). 


We finish the treatment of the case of $(1,4)$-affine 
control systems by writing the local normal form of such 
system, up to state-feedback transformation, in terms of 
the tuple of their primary invariants: Let $N$ be the 
solution of the following non-homogeneous second order 
linear ordinary differential equation 
w.r.t. the variable $x_1$ with prescribed initial values   
\begin{equation} 
\label{nonhom}\left\{  
\begin{aligned}
~&\frac{\partial^2 N}{\partial x_1^2}+{\mathcal 
I}_2\frac{\partial N}{\partial x_1}-{\mathcal I}_1 N-1=0;\\ 
~&N(0,x_2,x_3,x_4)\equiv 0,\quad \frac{\partial N}{\partial 
x_1}(x_1,x_2,x_3,x_4)\Bigl |_{x_1=0}\Bigr. \equiv 0, 
\end{aligned}\right.
\end{equation}
and the functions $\rho_1$, $\rho_2$ be the solution of the 
following homogeneous second order linear ordinary 
differential equations 
w.r.t. the variable $x_1$ with prescribed initial values 
\begin{equation} 
\label{hom}\left\{ 
\begin{aligned}
~&\frac{\partial^2 \rho_i}{\partial x_1^2}+{\mathcal 
I}_2\frac{\partial \rho_i}{\partial x_1}-{\mathcal I}_1 
\rho_i=0, \quad i=1,2,\\ 
~&\begin{pmatrix}\rho_1(0,x_2,x_3,x_4)&\rho_1(0,x_2,x_3,x_4)\\ 
\frac{\partial}{\partial x_1}\rho_1(0,x_2,x_3,x_4)& 
\frac{\partial}{\partial x_1}\rho_2(0,x_2,x_3,x_4) 
\end{pmatrix}\equiv\begin{pmatrix}1&0\\0&1\end{pmatrix}
\end{aligned}\right..
\end{equation}
Let also
\begin{equation}
\label{BK} B_k(x_2,x_3,x_4)=-\delta_{3k}+x_2 
\left(-\delta_{4k}+\int_0^{x_3}\psi_k(y,x_4)\,dy\right)+
\int_0^{x_2} (x_2-y)\beta_k(y,x_3,x_4)\, dy 
\end{equation}
for $1\leq k\leq 4$ (actually 
$B_k(x_2,x_3,x_4)=\frac{\partial a_k}{\partial 
x_1}(0,x_2,x_3,x_4)$, where the functions $a_k$ are as in 
(\ref{nf4})). Then a four-dimensional control-affine system 
(\ref{affcs1}) with the tuple of the primary invariants 
$({\mathcal I}_1, {\mathcal 
I}_2,\beta_1,\beta_2,\beta_3,\beta_4,\psi_1,\psi_2,\psi_3,\psi_4)$ 
at the point $q_0$ is state-feedback equivalent to the 
following system: 
\begin{equation} \label{rhoN} \left\{ 
\begin{aligned}
~&\dot x_1=1+(N+B_1\rho_2)u\\ ~&\dot 
x_2=(\rho_1+B_2\rho_2)u\\ ~&\dot x_i=B_i\rho_2 u ,\quad 
i=3,4 
\end{aligned}\right.,\qquad u\in{\mathbb R}. 
\end{equation} 

{\bf b) The case $n\geq 5$.} By (\ref{Fn}) and (\ref{N1}), 
in the canonical coordinates the vector fields $F_0$ and 
$F_1$ have the following form: \begin{equation} \label{nfn} 
 F_0=\sum_{k=1}^n 
a_k \frac{\partial}{\partial x_k}, \quad 
F_1=\frac{\partial}{\partial x_1}, \end{equation} where the 
components $a_k$ of $F_0$ satisfy the following system of 
partial differential equations 
\begin{equation}
\label{condnfn} 
\begin{split} &{\mathcal I}_{n-2}\frac{\partial^{n-2}a_k}{\partial 
x_1^{n-2}}+\sum_{j=3}^{n-3}{\mathcal I}_j \frac{\partial^j 
a_k}{\partial x_1^j}+\sum_{l=1}^n\Bigl(a_l \frac{\partial^2 
a_k}{\partial x_1^2}-\frac{\partial a_l}{\partial 
x_1}\frac{\partial a_k}{\partial x_l}\Bigr) +{\mathcal 
I}_2\frac{\partial a_k}{\partial x_1}+a_k+{\mathcal 
I}_1\delta_{1,k}=0, 
\\
&k=1,\ldots, n,  
\end{split}
\end{equation}
with the following restrictions on the boundary conditions 
for any $1\leq k\leq n$ 
\begin{subequations}
\label{boundn}  
\begin{gather}
  a_k(0,x_2,\ldots, x_n)\equiv \delta_{2k},\label{boundn:a}\\ 
 \cfrac{\partial a_k}{\partial 
x_1}(0,0,x_3\ldots,x_n)
\equiv \delta_{3k}, \label{boundn:b}\\ \cfrac{\partial^2 
a_k}{\partial x_1\partial x_2}(0,0,0,x_4,\ldots ,x_n)
\equiv -\delta_{4k}, \label{boundn:c}\\ \cfrac{\partial^j 
a_k}{\partial x_1^j}(0,\ldots,0,x_{j+1},\ldots,x_n) 
\equiv \delta_{j+3,k},\,\,2\leq j\leq n-3, \label{boundn:d} 
\end{gather}  
 \end{subequations} 
where $\delta_{ij}$ is the Kronecker symbol. Note also that 
the genericity assumption (\ref{gen3}) implies that
\begin{equation}
\label{lead} {\mathcal I}_{n-2}\neq 0. \end{equation} Let 
us introduce the following functions for any $1\leq k\leq 
n$:
\begin{subequations}
\label{primen} \begin{align} 
~&\beta_k(x_2,\ldots,x_n)\stackrel{def}{=}\cfrac{\partial^3 
a_k}{\partial x_1 \partial 
x_2^2}(0,x_2,\ldots,x_n),\label{primen:a}\\ 
~&\psi_k(x_3,\ldots,x_n)\stackrel{def}{=} \cfrac{\partial^3 
a_k}{\partial x_1 
\partial x_2\partial x_3}(0,0,x_3,\ldots,x_n)\label{primen:b},
\\ 
~&\phi_{kjl}(x_l,\ldots,x_n)\stackrel{def}{=}\cfrac{\partial^{j+1} 
a_k}{\partial x_1^j \partial x_l}a_k(0,\ldots,0,x_l,\ldots 
x_n),\quad 2\leq j\leq n-3,\,\,2\leq l\leq 
j+2\label{primen:c}. 
\end{align}
\end{subequations}  
So, with any germ at $q_0$ of an $n$-dimensional affine 
system (\ref{affcs1}), satisfying genericity assumptions 
(\ref{gen1}) and (\ref{gen2}), one can associate the 
ordered tuple 
\begin{equation}
\label{princen} \begin{split}\Big(\{{\mathcal 
I}_s(x_1,\ldots,x_n) 
\}_{s=1}^{n-2},\{\beta_{k}(x_2,\ldots,x_n)\}_{k=1}^n,
\{\psi_{k}(x_3,\ldots,x_n)\}_{k=1}^n,\\ 
\{\phi_{kjl}(x_l,\ldots, x_n):1\leq k\leq n, 2\leq j\leq 
n-3, 2\leq l\leq j+2\} \Big) 
\end{split}
\end{equation}
of state-feedback invariants.
of functions of $n-1$ variables at $0$. We call it {\it the 
tuple of the primary invariants of the $(1,n)$-affine 
control system} (\ref{affcs1}) {\it with $n>4$ at the point 
$q_0$}. Note that by (\ref{prime4}) for any $1\leq k\leq n$ 
the functional invariants $\beta_k$ and $\psi_k$ have the 
weight $3$ , the functional invariants $\phi_{kjl}$ have 
the weight $j+1$, while by \eqref{condnf4} and 
\eqref{bound4} for any $1\leq s\leq n-2$ the functional 
invariants ${\mathcal I}_s$ have the weight $n-2$. 

Further, fixing $\beta_k$ and $\psi_k$ and using 
\eqref{boundn:b} and \eqref{boundn:c}, one can find 
$\frac{\partial a_k}{\partial x_1}(0,x_2,\ldots,x_n)$ for 
any $1\leq k\leq n$ by the appropriate integrations. 
Similarly, fixing $\{\phi_{kjl}\}_{l=2}^{j+2}$ for given 
$j$, $2\leq j\leq n-3$, and using \eqref{boundn:d}, one can 
find $\frac{\partial^j a_k}{\partial 
x_1^j}(0,x_2,\ldots,x_n)$ for any $1\leq k\leq n$ by the 
appropriate integrations. Finally, if we suppose that all 
functions $\beta_k$, $\psi_k$, and $\phi_{kjl}$ are real 
analytic and fix also real analytic $\{{\mathcal 
I}_s\}_{s=1}^{n-2}$, then from the knowledge of 
$\frac{\partial^j a_k}{\partial x_1^j}(0,x_2,\ldots,x_n)$ 
for all $1\leq j\leq n-3$, condition \eqref{boundn:a}, and 
differential equation (\ref{condnfn}) we can recover the 
functions $a_k(x_1,\ldots,x_n)$ and therefore our affine 
control system itself, just using the classical Cauchy - 
Kowalewsky theorem for system (\ref{condnfn}). We summarize 
all above in the following: 


\begin{theor}
\label{classn} 
If $n\geq 5$, then given an arbitrary tuple \eqref{princen} 
of real analytic functions $n-2$ germs $\{{\mathcal 
I}_j\}_{j=1}^{n-2}$ of real analytic functions there exists 
a unique, up to state-feedback real analytic transformation 
of the type (\ref{feedback}), $n$-dimensional real analytic 
control-affine system with scalar input, satisfying genericity 
assumptions (\ref{gen1}),(\ref{gen2}), and 
(\ref{gen3}),such that the tuple \eqref{princen} is its 
tuple of the primary invariants at the given point $q_0$.
 In other words, the 
tuples of the primary invariants give the regular 
$(n-2,2)$-parameterization of the considered classification 
problem in the real analytic category with the $(n-1)\times 
(n-4)$ parameterization matrix $P$ such that for $n=5$ 
\begin{equation}
\label{parmat5} P=(5,10,10,3)^T \end{equation} and for 
$n>5$ 
\begin{equation}
\label{parmatn} P=\begin{pmatrix} 0&\hdotsfor{3}&0&n~~~~~\\ 
\vdots&~& ~&\sddots&\sddots&\vdots~~~~~\\ 
\vdots&~&\sddots&\sddots&~&\vdots~~~~~\\ 
\vdots&\sddots&\sddots&~&~&\vdots~~~~~\\ 
0&n&\hdotsfor{3}&n~~~~~\\ n&n& \hdotsfor{3}&n~~~~~\\2n&n& 
\hdotsfor{3}&n~~~~~\\2n&n& 
\hdotsfor{3}&n~~~~~\\0&\hdotsfor{3}&0&n-2 
\end{pmatrix}
\end{equation} 
(in the matrix $P$ all entries in the triangle with 
vertices in the $(1,1)$-entry, the $(1,n-5)$-entry, and the 
$(n-5,1)$-entry are equal to $0$, while all entries in the 
triangle with vertices in the $(1,n-4)$-entry, the 
$(n-4,1)$-entry, and the $(n-4,n-4)$-entry are equal to 
$n$). The Poincare series $M(t)$ of the considered 
classification problem satisfies 
\begin{equation}
\label{poinn} 
\begin{split}
M(t)=nt^3\left(\cfrac{2}{(1-t)^n}+\cfrac{2}{(1-t)^{n-1}}+
\cfrac{1}{(1-t)^{n-2}}\right)\\ 
+n\sum_{i=4}^{n-2}t^i\sum_{j=n-i}^{n-1}\cfrac{1}{(1-t)^{j+1}}
+\cfrac{(n-2)t^{n-2}}{(1-t)^{n+1}}.
\end{split}
\end{equation} 
\end{theor} 

It turns out that $(n-2,2)$ is the characteristic pair of 
the considered classification problem. Indeed, let as 
before $w_0$ be the order of zero of $M(t)$ at $t=0$, $N$ 
is the order of pole of $(1-t)M(t)$ at $t=1$, and $d$ is 
the degree of $M(t)$ (at infinity). Then from (\ref{poin4}) 
it follows that $w_0=3$, $N=n$, and $d=n-5$. Hence from 
(\ref{Nicein3}) it follows that 
\begin{equation}\label{NSn1}
{\rm NS}\bigl(M(t)\bigr)\subset\{(w,l):w\geq n-3\} . 
\end{equation}
Further, using \eqref{Nicein2}, one has 
\begin{equation}\label{NSn2}
\begin{split}
&{\rm NS}\bigl(M(t)\bigr)\cap\{(w,l):w= 
2\}\subset\{(n-3,1)\},\\&{\rm 
NS}\bigl(M(t)\bigr)\cap\{(w,l):w=n-2\}\subset\{(n-2,1), 
(n-2,2)\} . 
\end{split}
\end{equation} 
But by the previous theorem our classification problem 
admits $(n-2,2)$-parameterization such that its 
parameterization metric has a positive entry in the 
right-upper corner. Therefore from Corollary \ref{usecor} 
it follows that $(n-3,1)\notin {\rm NS}\bigl(M(t)\bigr)$. 
Since $(n-2,2)\prec(n-2,1)$ we can conclude from 
(\ref{NSn1}) and (\ref{NSn2}) that $(n-2,2)$ is the minimal 
element of ${\rm NS}\bigl(M(t)\bigr)$. In other words, 
$(n-2,2)$ is the characteristic pair of our classification 
problem. Also, it implies that the characteristic matrix 
${\mathcal C}$ of the problem is equal to ${\rm Norm}(P)$, 
where $P$ is as in Theorem \ref{classn}. If $n=5$ then 
obviously ${\mathcal C}={\rm Norm}(P)=P$. For $n>5$ one can 
calculate all nontrivial entries of ${\rm Norm}(P)$, using 
identities (\ref{norm}). It gives all nontrivial entries of 
the characteristic matrix ${\mathcal C}$: 
\begin{equation}
\label{charn}
\begin{split}
&{\mathcal C}_{n-1, n-4}=n-2,\\ &{\mathcal 
C}_{n-2,n-4}=n(n-3),\\& {\mathcal C}_{i, 
n-4}=n\Bigl(1+\sum_{l=1}^{n-3-i}\tbinom{i+2l-3}{l}+\tbinom{2n-9-i}{n-3-i} 
+2\tbinom{2n-8-i}{n-2-i}+\tbinom{2n-7-i}{n-7-i}\Bigr),\quad 
2\leq i\leq n-3,\\ &{\mathcal C}_{1,n-4}=n,\\ &{\mathcal 
C}_{1j}=n\Bigl(\tbinom{n-4+j}{j-1}+ 
2\tbinom{n-5+j}{j-1}+2\tbinom{n-6+j}{j-1}+
\tbinom{n-7+j}{j-1}-1-\\&-\sum_{l=1}^{j-1}\tbinom{n+2l-7-j}{l}\Bigr),\quad 
1\leq j\leq n-5. 
\end{split}
\end{equation} 
Recall that ${\mathcal C}_{ij}$ is the intrinsic number of 
the functional invariants of $i+1$ variables and the weight 
$j+2$. In order to obtain a characteristic parameterization 
from the parameterization by the tuple of the primary 
invariants one can implement some series of rearrangement 
of the primary invariants according to the series of 
elementary transformations from the matrix $P$ to ${\rm 
Norm}(P)$, as was described in the previous section (see, 
for example, the formula (\ref{FG}) and the paragraph after 
it). 

\begin{remark} \label{affrem} {\rm  
Two control systems $\dot y={\mathcal F}(y,v)$ and $\dot 
{\tilde y}=\widetilde {\mathcal F}(\tilde y,\tilde v)$, 
with $m$-dimensional state-space $S$ and one dimensional 
control space $V$ are called {\it micro-locally 
state-feedback equivalent at the point $(y_0, v_0)\in 
S\times V$}, if there exist the state-feedback 
transformation 
\begin{equation*}
\label{feedbackn} \left\{\begin{array} {l} \tilde 
y=\Phi_1(y)\\ \tilde v= \Phi_2(y,v)\\ 
y_0=\Phi_1(y_0),\,\,v_0=\Phi_2(y_0,v_0)\end{array}\right. 
\end{equation*} 
such that in the neighborhood of $(y_0,v_0)$ in $S\times V$ 
the following identity holds: $$\widetilde {\mathcal 
F}\bigl(\Phi_1(y),\Phi_2(y,v)\bigr)=d \Phi_1{\mathcal 
F}(y,v).$$ The affine $(1,m+1)$ control system 
(\ref{affcor}) will be called {\it the affine extension} of 
the control system (\ref{nonaff}). It is not difficult to 
show that two control systems with scalar input are 
micro-locally state-feedback equivalent at some point 
$(y_0, v_0)\in S\times V$ if and only if their affine 
extensions are locally equivalent w.r.t. the state-feedback 
transformations of the type (\ref{feedback}) at the same 
point. Note also that the affine extensions of generic 
$m$-dimensional control systems with scalar input 
are generic in the set of all $(1,m+1)$-affine systems. 
Using this fact and Theorems \ref{class4}, \ref{classn} one 
obtains the micro-local parameterization of non-affine 
$m$-dimension control systems with scalar inputs and $m\geq 3$ 
by the tuples of the primary invariants of their affine 
extensions (in $C^\infty$ category for $m=3$ and $C^\omega$ 
category for $m\geq 4$). Obviously, the Poincare series, 
the characteristic pair and the characteristic matrix of 
the micro-local state-feedback classification problem for 
$m$-dimensional control systems with scalar input 
are exactly the same as in the case of the state-feedback 
classification problem of $(1,m+1)$ control-affine systems. 
Besides, since in our method of normalization of 
$(1,n)$ control-affine system with $n\geq 5$ we rectify the 
vector field $f_1$, then in the case $m\geq 4$ a generic 
$m$-dimensional real analytic control system with 
the prescribed tuple (\ref{princen}) of the primary 
invariants of its affine extension is micro-locally 
state-feedback equivalent in the real analytic category to 
the system $$ \dot y_s=f_s(y_1,\ldots,y_m, v),\quad 1\leq 
s\leq m, $$ such that $f_s(x_2,x_3,\ldots x_m, 
x_1)=a_{s+1}(x_1,\ldots x_{m+1})$, where the tuple 
$\{a_k(x_1,\ldots,x_m)\}_{k=1}^{m+1}$ is the solution of 
the system of partial differential equations 
(\ref{condnfn}) with boundary conditions, which can be 
expressed by the primary invariants $\beta_k$, $\psi_k$, 
and $\phi_{kjl}$, using (\ref{primen}).} 
\end{remark}
 
\section
{Reduction of control-affine systems with two-dimensional input to the 
scalar input case in dimensions four and five} \indent 
\setcounter{equation}{0} 

For $r=2$ the system (\ref{affcs}) has the following form
\begin{equation}
\label{affcs2} \dot q=f_0(q)+u_1\,f_1(q)+u_2 \,f_2(q), 
\quad q\in M,\quad u_1,\,u_2\in\mathbb R.\end{equation} Our 
aim is to assign to the system (\ref{affcs2}) in a 
canonical way an affine subsystem with scalar input\footnote 
 {the meaning of the word 
``subsystem' is that at any point $q$ the set of its 
admissible velocities is a subset of the set of the 
admissible velocities of the original system at $q$.}. It 
turns out that in the case $n=4$ the original system can be 
recovered from it uniquely up to a feedback transformation, 
while in the case $n=5$ such unique recovering is possible 
after introducing an additional invariant function of $n$ 
variables (which is natural in view of the estimates for 
the number of functional parameters, given in the 
Introduction). 

{\bf 4.1 Preliminaries.} Let us look on (\ref{affcs2}), as 
on the time optimal control problem and find its extremals. 
First we introduce some notations. Let $T^*M$ be the 
cotangent bundle of $M$ with canonical symplectic form 
$\sigma$. Denote by $h_i$, $0\leq i\leq 2$, the following 
functions on $T^*M$: 
\begin{equation}
\label{quasi25} h_i(\lambda)=p\cdot 
f_i(q),\,\,\lambda=(p,q),\,\, q\in M,\,\, p\in T_q^* M.
\end{equation}
For a given function $G:T^*M\mapsto \mathbb R$ denote by 
$\vec G$ the corresponding Hamiltonian vector field defined 
by the relation $\sigma(\vec G,\cdot)=-d\,G(\cdot)$. For a 
given vector distribution $D$ on $M$ (i.e., a subbundle of 
the tangent bundle), define the $l$th power $D^l$ by the 
recursive relation 
 $$D^l=D^{l-1}+[D,D^{l-1}], \quad D^1=D,$$ and denote by $(D^l)^\perp\subset T^*M$ 
the annihilator of $D^l$, namely $$(D^l)^\perp=\{(p,q)\in 
T^*M: p\cdot v=0 \,\,\forall v\in D^l(q)\}.$$ 

In 
the introduced notations the Hamiltonian of Pontryagin 
Maximum Principle for the time optimal problem 
(\ref{affcs2}) can be written as follows: 
\begin{equation}
\label{Hamu}
 H(\lambda, u_1,u_2)= h_0(\lambda)+u_1 h_1(\lambda)+u_2 h_2(\lambda),\quad 
 \lambda\in T^*M,\,\, 
 u_1,u_2\in{\mathbb R}.
 \end{equation}
Let $\gamma(\cdot)$ be an extremal of (\ref{affcs2}) with 
extremal control functions $\bar u_1(t)$ and $\bar u_2(t)$. 
Then 
\begin{equation}
\label{dotl} 
 \dot\gamma(t)= \vec h_0\bigl(\gamma(t)\bigr)+
 \bar u_1(t) \vec h_1\bigl(\gamma(t)\bigr)+\bar u_2(t) \vec 
 h_2\bigl(\gamma(t)\bigr)
\end{equation} 
and from the maximality condition for $H$ it follows that 
\begin{equation}
\label{belong} \gamma(\cdot) \subset \{\lambda\in T^*M: 
h_1(\lambda)=h_2(\lambda)=0\}. 
\end{equation}
If we denote $D_2={\rm span}(f_1,f_2)$, then (\ref{belong}) 
is equivalent to $\gamma(\cdot)\subset (D_2)^\perp$. 
Combining (\ref{dotl}) and (\ref{belong}), we obtain 
\begin{equation}
\label{deriv}
d_{\gamma(t)}h_i\bigl(\dot\gamma(t)\bigr)=0, i=1,2.
\end{equation} 

Then from  (\ref{dotl}) and (\ref{deriv}) it follows  
\begin{equation}
\label{pois} 
 \begin{split}
 &
\{h_0,h_1\}\bigl(\gamma(t)\bigr)+
\bar u_2(t)\{h_2,h_1\}\bigl(\gamma(t)\bigr)=0,\\ 
&\{h_0,h_2\}\bigl(\gamma(t)\bigr)+\bar u_1(t)\{h_1,h_2\}\bigl(\gamma(t)\bigr)=0
\end{split} 
\end{equation}
(here $\{h_i,h_j\}$ are Poisson brackets of the 
Hamiltonians $h_i$ and $h_j$: $\{h_i,h_j\}=dh_j(\vec 
h_i)$). 
Now suppose  that \begin{equation} \label{denote} 
\dim\,D_2^2=3
\end{equation}
Then relations (\ref{pois}) implies that the extremals of 
(\ref{affcs2}), lying in $(D_2)^\perp\backslash 
(D_2^2)^\perp$, are exactly the integral curves of the 
vector field 
\begin{equation}
\label{hamaf} \vec X=\vec 
h_0+\cfrac{\{h_0,h_2\}}{\{h_2,h_1\}} \vec 
h_1+\cfrac{\{h_1,h_0\}} {\{h_2,h_1\}}\vec h_2 
\end{equation}
(which is the Hamiltonian vector field, corresponding to 
the Hamiltonian $X=h_0+\cfrac{\{h_0,h_2\}}{\{h_2,h_1\}} 
h_1+\cfrac{\{h_1,h_0\}} {\{h_2,h_1\}}h_2$).
Denote by $V$ the affine subbundle of $TM$, defined by 
system (\ref{affcs2}) and $V(q)$ be the set of all 
admissible velocities of the system (\ref{affcs2}) at the 
point $q$, 
$$V(q)=\{f_0(q)+u_1 f_1(q)+u_2f_2(q): u_1,u_2\in{\mathbb 
R}\}.$$ 
Let $\pi:T^*M\mapsto M$ be the canonical projection. The 
set \begin{equation} \label{ext} {\rm Ext}(q)=\{\pi_*\vec 
X(\lambda):\lambda\in T_q^*M\cap (D_2)^\perp\backslash 
(D_2^2)^\perp\}, \quad q\in M 
\end{equation}
is the subset of $V(q)$,
consisting of the velocities of all extremal trajectories 
starting at $q$ and having a lift in $(D_2)^\perp\backslash 
(D_2^2)^\perp$. 

Among all extremals on $(D_2)^\perp\backslash (D_2^2)^\perp$, one can distinguish  so-called abnormal extremals, i.e., the extremals lying on the zero level set of the Hamiltonian $X$. 
Denote $D_3={\rm span}(f_0,f_1,f_2)$ and suppose that
\begin{equation}
\label{D23} \dim \,(D_2^2+D_3)=4 
\end{equation}
The set 
\begin{equation} 
\label{extab} {\rm Abn}(q)=\{\pi_*\vec 
X(\lambda):\lambda\in T_q^*M\cap (D_3)^\perp\backslash 
(D_2^2)^\perp\}, \quad q\in M 
\end{equation}
is the subset of ${\rm Ext}(q)$,
consisting of the velocities of all abnormal extremal 
trajectories starting at $q$ and having a lift in 
$(D_2)^\perp\backslash (D_2^2)^\perp$. One can show that 
for generic affine systems of the type (\ref{affcs2}) ${\rm 
Ext}(q)=V(q)$ in the case $n\geq 5$ and ${\rm Abn}(q)=V(q)$ 
in the case $n\geq 6$. But in the case $n=4$ and $n=5$ 
either ${\rm Ext}(q)$ or ${\rm Abn}(q)$ (or both of them) 
define the proper subsystem of the original system 
(\ref{affcs2}). Moreover, it turns out that these 
subsystems are affine with scalar input, so one can apply the 
theory of the previous section. Now let us consider the 
cases $n=4$ and $n=5$ separately. 

{\bf 4.2 The case $n=4$.} Let 
$$[V,D_2](q)=\{[X,Y](q):\,X\in V,\,Y\in D_2,\,\,{\rm 
are}\,\, {\rm vector}\,\,{\rm fields}\}.$$ It is not 
difficult to show that $[V,D_2](q)$ is a linear space and 
\begin{equation}
\label{VD2} [V,D_2](q)={\rm 
span}\bigl(f_1(q),f_2(q),[f_1,f_2](q),[f_0,f_1](q),[f_0,f_2](q)\bigr)
\end{equation}
The crucial observation is formulated in the following 

\begin{prop}
\label{afflem} The set ${\rm Ext}(q)$ is an affine line, 
provided that (\ref{denote}) holds and  
\begin{equation}
\label{condred} \dim\, [V,D_2](q)=4
\end{equation} 
\end{prop}

{\bf Proof.} Take some vector field $f_3$ such that the 
tuple $(f_0,f_1,f_2,f_3)$ constitutes the frame on $M$. 
Denote by $c_{ji}^k$ the structural functions of this 
frame, i.e., the functions, satisfying 
\begin{equation}
\label{struct} [f_i, f_j]=\sum_{k=0}^3c_{ji}^kf_k. 
\end{equation}
Using the following well-known property of the Poison 
brackets 
\begin{equation}
\label{poisdef} \{h_i,h_j\}(p,q)=p\cdot[f_i,f_j](q),\quad 
q\in M, p\in T^*_qM \end{equation} 
and (\ref{hamaf}), one can easily obtain that 
\begin{equation}
\label{ext1} {\rm Ext}(q)= \Pi(q)\cap V(q), \end{equation} 
where 
\begin{equation}
\label{pi} \begin{split}&\Pi(q)=\{\bigl(c_{12}^0(q)\nu 
+c_{12}^3(q) \mu\bigr) f_0(q)+\bigl(c_{20}^0(q)\nu 
+\bigr.\\ &~\quad~ \,\,\,\,\,\, \bigl.c_{20}^3(q)\mu\bigr) 
f_1(q) + \bigl(c_{01}^0(q)\nu+ 
c_{01}^3(q)\mu\bigr)f_2(q):\mu,\nu\in \mathbb R\} 
\end{split} 
\end{equation}
From assumption (\ref{condred}) and identity (\ref{VD2}) it 
follows that 
$\Pi(q)$ is a plane. Assumption (\ref{denote}) implies that 
the plane $\Pi(q)$ is not parallel to the plane $V(q)$. 
Note also both $\Pi(q)$ and $V(q)$ belong to $D_3(q)$. 
Hence by (\ref{ext1}) the set ${\rm Ext}(q)$ is an affine 
line. 
$\Box$ \vskip .2in 

Consider the control system such that ${\rm Ext}(q)$ is its 
set of the admissible velocities at $q$. By Proposition 
\ref{afflem} it is an affine system with scalar input. We call 
this system {\it the reduction of the four-dimensional 
control-affine system} (\ref{affcs2}). The following 
proposition gives another characterization of the reduction 
of the system (\ref{affcs2}): 

\begin{prop}\label{charlemma}
Assume that the four-dimensional control-affine system 
(\ref{affcs2}) satisfies the conditions (\ref{denote}) and 
(\ref{condred}). Then the subsystem 
\begin{equation}
\label{affcs21} \dot q=g_0+u g_1, \end{equation} of 
(\ref{affcs2}) is its reduction if and only if 
\begin{equation}
\label{chareq} [g_0,g_1]\in D_2. 
\end{equation}
\end{prop}

{\bf Proof.} By definition, the system (\ref{affcs21}) is 
the reduction of (\ref{affcs2}) if and only if 
\begin{equation}
\label{ext2} {\rm Ext}(q)=\{g_0(q)+t g_1(q): t\in \mathbb 
R\}. 
\end{equation} On the other hand, on can take from the 
beginning $f_0=g_0$ and $f_1=g_1$. Then comparing 
(\ref{ext2}) with (\ref{ext1}) and (\ref{pi}) we obtain 
that the system (\ref{affcs21}) is the reduction of 
(\ref{affcs2}) if and only if $c_{01}^0=c_{01}^3=0$, which 
is equivalent to $[g_0,g_1]\in {\rm span}(g_1,f_2)=D_2$. 
$\Box$ 

\begin{cor}
\label{recover} Assume that a four-dimensional affine 
control system (\ref{affcs21}) satisfies 
\begin{equation} \label{geng1}
 \dim\,{\rm 
 span}\bigl(g_1,[g_1,g_0],\bigl[g_1,[g_1,g_0]\bigl],
 \bigl[g_0,[g_1,g_0]\bigr]\bigr)=4
 \end{equation}
 Then any  
 four-dimensional control-affine system with two-dimensional input, 
 having the system (\ref{affcs21}) as its reduction, is 
 feedback equivalent to the system
 \begin{equation}
\label{affcsrec} \dot q=g_0+u_1 g_1+u_2[g_0,g_1],\quad 
u_1,u_2\in \mathbb R. 
\end{equation} 
 \end{cor}
 
 {\bf Proof.} 
 First by assumption (\ref{geng1}) and relation (\ref{VD2}) 
 the system 
 (\ref{affcsrec}) satisfies conditions (\ref{denote}) and 
 (\ref{condred}) (where $f_0$, $f_1$, and $f_2$ are replaced by 
 $g_0$, $g_1$ and $[g_0,g_1]$). Hence, by Proposition \ref{afflem} 
 the system (\ref{affcsrec}) 
 admits the reduction and by Proposition \ref{charlemma} this reduction is the system 
 (\ref{affcs21}). 
 On the other hand, suppose that some system  
 (\ref{affcs2}) has the reduction  (\ref{affcs21}). Then from the previous proposition 
 $[g_0,g_1]\in {\rm span}( g_1, f_2)$. According to  
 (\ref{geng1}), $g_1$ and $[g_0,g_1]$ are linearly 
 independent. Hence the system (\ref{affcs2}) is feedback 
 equivalent to (\ref{affcsrec}). $\Box$
 \vskip .2in

According to the previous proposition, a generic 
four-dimensional control-affine system with two-dimensional input can 
be uniquely, up to a feedback transformation, recovered 
from its reduction. Suppose that the reduction 
(\ref{affcs21}) of the system (\ref{affcs2}) satisfies 
(\ref{geng1}) and 
\begin{equation} \label{geng2}
 \dim\,{\rm 
 span}\bigl(g_0,g_1,[g_1,g_0],\bigl[g_1,[g_1,g_0]\bigl]
 \bigl)=4.
 \end{equation} 
Then we can apply to the system (\ref{affcs21}) all 
constructions of section 2. In particular, one can 
construct the tuple of the primary invariants of 
(\ref{affcs21}) at a given point, which are also feedback 
invariants of the original system (\ref{affcs2}). Note that 
the set of germs of systems of the type (\ref{affcs2}) 
having the reductions, which satisfies conditions 
(\ref{geng1}) and (\ref{geng2}), is generic. Combining 
Theorem \ref{class4}, Corollary \ref{recover}, and normal 
form (\ref{rhoN}), we obtain the following classification 
of generic germs of systems of the type (\ref{affcs2}) in 
terms of the tuple of the primary invariants of their 
reductions: 
\begin{theor}
\label{class24
} 
Given two arbitrary germs ${\mathcal I}_1$ and ${\mathcal 
I}_2$ of functions of four variables at $0$, four arbitrary 
germs $\beta_1,\beta_2,\beta_3,\beta_4$ of functions of 
three variables at $0$, and four arbitrary germs 
$\psi_1,\psi_2,\psi_3,\psi_4$ of functions of two variables 
at $0$ there exists a unique, up to state-feedback 
transformation of the type (\ref{feedback}), 
four-dimensional control-affine system with two-dimensional input such 
that its reduction 
satisfies genericity assumptions (\ref{geng1}) and 
(\ref{geng2}) and $({\mathcal I}_1, {\mathcal 
I}_2,\beta_1,\beta_2,\beta_3,\beta_4,\psi_1,\\\psi_2,\psi_3,\psi_4)$ 
is the tuple of the primary invariants of the reduction at 
the given point $q_0$. This control system is 
state-feedback equivalent to the following one: 
\begin{equation} \label{rhoN2} \left\{ 
\begin{aligned}
~&\dot x_1=1+(N+B_1\rho_2)u_1+\Bigl(\frac{\partial 
N}{\partial x_1} +\beta_1\frac{\partial\rho_2 }{\partial 
x_1}\Bigr)u_2\\ ~&\dot 
x_2=(\rho_1+B_2\rho_2)u_1+\Bigl(\frac{\partial\rho_1}{\partial 
x_1}+B_2\frac{\partial\rho_2 }{\partial x_1}\Bigr)u_2\\ 
~&\dot x_i=B_i\rho_2 u_1+\beta_i\frac{\partial\rho_2 
}{\partial x_1} u_2 ,\quad i=3,4 
\end{aligned}\right.,\quad u_1,u_2\in{\mathbb R}, 
\end{equation}
where $N$ is the solution of (\ref{nonhom}), $\rho_i$, 
$i=1,2$ are the solutions of (\ref{hom}), and $B_k$, $1\leq 
k\leq 4$, are as in (\ref{BK}). The Poincare series, the 
characteristic pair and the characteristic matrix of the 
classification problem are exactly the same as in the case 
of $(1,4)$ control-affine systems. 
\end{theor}

\begin{remark}
{\rm It is easy to show that in the case $n=4$ the set 
${\rm Abn}(q)$ consists of one vector provided that 
(\ref{D23}) holds. Besides, if the system (\ref{affcs21}) 
is the reduction of the system (\ref{affcs2}) and it 
satisfies (\ref{geng1}) and (\ref{geng2}), then ${\rm 
Abn}(q)$ is exactly its canonical drift.}$\Box$ 
\end{remark} 

\begin{remark}
\label{way2} {\rm Actually, there is another intrinsic way 
to assign to the system (\ref{affcs2}), satisfying 
(\ref{D23}), an affine subsystem with scalar input: As a drift 
one can take again $Ab(q)$. It remains to define 
canonically the direction of the affine line of the 
reduction. For this note first that the distribution $D_2$ 
satisfies (\ref{denote}) because of assumption (\ref{D23}). 
Therefore through any point of $M$ the unique 
(unparameterized) abnormal extremal trajectory of the rank 2 
distribution $D_2$ passes: the line subdistribution $L$ of 
$D_2$, tangent to the abnormal extremal trajectories at any 
point is characterized by the relation $[L,D^2]\subseteq 
D^2$. The direction of the affine line of the reduction can 
be taken parallel to $L$. 
The direction of $L$ is different in general from the 
direction of the affine line in the first reduction. But 
this new reduction is worse than the previous one, because 
the original system (\ref{affcs2}) is not uniquely 
recovered from it: if $(\bar g_0,\bar g_1)$ is the 
canonical pair of the new reduction (by construction and 
the previous remark $g_0={\rm Abn}$), then the field $f_2$ 
can be taken in the form $f_2=\alpha \bar g_0 +[\bar 
g_1,\bar g_0]$, where the function $\alpha$ satisfies some 
second order ordinary differential equation along each 
integral curve of $\bar g_1$. 
Note that 
the direction $\bar g_1$ depends on the second jet of the 
original system (\ref{affcs2}), while the direction of the 
affine line in the first reduction depends only on the 
first jet. This could be the reason for the loss of some 
information about the original system during the reduction 
described in the present remark.} $\Box$ 
\end{remark} 

{\bf 4.3 The case $n=5$.} In this case by analogy with 
Proposition \ref{afflem} we have

\begin{prop}
\label{affab} The set ${\rm Abn}(q)$ is an affine line, 
provided that (\ref{D23}) holds and \begin{equation} 
\label{D3sqr} \dim\,D_3^2=5. \end{equation} 
\end{prop}

{\bf Proof.} Take some vector fields $f_3$ and $f_4$ such 
that the tuple $(f_0,f_1,f_2,f_3, f_4)$ constitutes the 
frame on $M$. 
By analogy with (\ref{struct}), let $c_{ji}^k$, $0\leq 
i,j,k\leq 4$, be the structural functions of this frame. 
From (\ref{hamaf}) and (\ref{extab}), using (\ref{poisdef}) 
and the fact that $h_0=0$ on $(D_3)^\perp$, one can easily 
obtain that 
\begin{equation}
\label{abn1} {\rm Abn}(q)= \Pi_1(q)\cap V(q), 
\end{equation} where 
\begin{equation}
\label{pi1} \begin{split}&\Pi_1(q)=\{\bigl(c_{12}^3(q)\nu 
+c_{12}^4(q) \mu\bigr) f_0(q)+\bigl(c_{20}^3(q)\nu 
+\bigr.\\ &~\quad~ \,\,\,\,\,\, \bigl.c_{20}^4(q)\mu\bigr) 
f_1(q) + \bigl(c_{01}^3(q)\nu+ 
c_{01}^4(q)\mu\bigr)f_2(q):\mu,\nu\in \mathbb R\} 
\end{split} 
\end{equation}
From assumption (\ref{D3sqr}) and identity (\ref{VD2}) it 
follows that 
$\Pi_1(q)$ is a plane. Assumption (\ref{D23}) implies that 
the plane $\Pi_1(q)$ is not parallel to the plane $V(q)$. 
Note also both $\Pi_1(q)$ and $V(q)$ belong to $D_3(q)$. 
Hence by (\ref{abn1}) the set ${\rm Abn}(q)$ is an affine 
line. $\Box$ \vskip .2in 

Consider the control system such that ${\rm Abn}(q)$ is its 
set of the admissible velocities at $q$. By Proposition 
\ref{affab} it is an affine system with scalar input. We call 
this system {\it the reduction of the five-dimensional 
control-affine system} (\ref{affcs2}). The following 
proposition gives another characterization of the reduction 
of the system (\ref{affcs2}): 

\begin{prop}\label{charabn}
Assume that the five-dimensional control-affine system 
(\ref{affcs2}) satisfies the conditions (\ref{D3sqr}) and 
(\ref{D23}). Then the subsystem 
\begin{equation}
\label{affcs25} \dot q=g_0+u g_1, \end{equation} of 
(\ref{affcs2}) is its reduction if and only if 
\begin{equation}
\label{chareq1} [g_0,g_1]\in D_3. 
\end{equation}
\end{prop}

{\bf Proof.} By definition, the system (\ref{affcs21}) is 
the reduction of (\ref{affcs2}) if and only if 
\begin{equation}
\label{abn2} {\rm Abn}(q)=\{g_0(q)+t g_1(q): t\in \mathbb 
R\}. 
\end{equation} On the other hand, on can take from the 
beginning $f_0=g_0$ and $f_1=g_1$. Then comparing 
(\ref{abn2}) with (\ref{abn1}) and (\ref{pi1}) we obtain 
that the system (\ref{affcs25}) is the reduction of 
(\ref{affcs2}) if and only if $c_{01}^3=c_{01}^4=0$, which 
is equivalent to $[g_0,g_1]\in {\rm span}(g_0, 
g_1,f_2)=D_3$. $\Box$ 

\begin{cor}
\label{recover5} Assume that a five-dimensional affine 
control system (\ref{affcs21}) satisfies 
\begin{equation} \label{geng15}
 \dim\,{\rm 
 span}\bigl(g_0,g_1,[g_1,g_0],\bigl[g_1,[g_1,g_0]\bigl],
 \bigl[g_0,[g_1,g_0]\bigr]\bigr)=5
 \end{equation}
 Then a 
 five-dimensional control-affine system with two-dimensional input, 
 has the system (\ref{affcs21}) as its reduction if and only if it is  
 feedback equivalent to the system
 \begin{equation}
\label{affcsrec5} \dot q=g_0+u_1 g_1+u_2(\alpha g_0+
[g_0,g_1]),\quad u_1,u_2\in \mathbb R, 
\end{equation} 
where $\alpha$ is some function.
\end{cor}

 {\bf Proof.} 
 First by assumption (\ref{geng15})  
 the system 
 (\ref{affcsrec5}) satisfies conditions (\ref{D3sqr}) and 
 (\ref{D23}). Hence, by Proposition \ref{affab} 
 the system (\ref{affcsrec5}) 
 admits the reduction and by Proposition \ref{charabn} 
 this reduction is the system 
 (\ref{affcs21}). 
 On the other hand, suppose that some system  
 (\ref{affcs2}) has the reduction  (\ref{affcs21}). Then from the previous proposition 
 $[g_0,g_1]\in {\rm span}(g_0, g_1, f_2)$. According to  
 (\ref{geng1}), $g_0$, $g_1$ and $[g_0,g_1]$ are linearly 
 independent. Hence, $f_2=\xi_0 g_0+\xi_1 g_1+\xi_2 
 [g_0,g_1]$, where $\xi_2\neq 0$. Hence by a feedback transformation we can 
 replace $f_2$ with $\alpha g_0+[g_0,g_1]$. $\Box$
 \vskip .2in
 
 So, in contrast to the case $n=4$,  
five-dimensional control-affine systems with two-dimensional input 
cannot be recovered from its reduction only. Suppose that 
the system (\ref{affcs2}) has the reduction (\ref{affcs25}) 
satisfying condition (\ref{geng15}) and also 
\begin{eqnarray} &~&\label{geng25}
 \dim\,{\rm 
 span}\bigl(g_0,g_1,[g_1,g_0],\bigl[g_1,[g_1,g_0]\bigl],
 \bigl[g_1,\bigl[g_1,[g_1,g_0]\bigl]\bigl]\bigl)=5,\\
 &~&\label{geng35}
 \dim\,{\rm 
 span}\bigl(g_1,\bigl[g_0,[g_0,g_1]\bigl],[g_1,g_0],\bigl[g_1,[g_1,g_0]\bigl],
 \bigl[g_1,\bigl[g_1,[g_1,g_0]\bigl]\bigl]\bigl)=5.
 \end{eqnarray} 
Then we can apply to the system (\ref{affcs25}) all 
constructions of section 2. In particular, let $(G_0,G_1)$ 
be the canonical pair of the system (\ref{affcs25}). As 
before, let $V$ be the affine subbundle of $TM$, defined by 
system (\ref{affcs2}). Then by the same arguments, as in 
the proof of Corollary \ref{recover5}, there exist a unique 
vector field $G_2\in V$ such that 
\begin{equation}
\label{recinv} G_2=R G_0+[G_0,G_1]. \end{equation} By 
construction, the function $R$ is a feedback invariant of 
system (\ref{affcs2}). Moreover, the system (\ref{affcs2}) 
can be uniquely, up to a feedback transformation, recovered 
from its reduction and the function $R$. 

Fix some point $q_0$ in $M$. Let $\Phi_5$ be as in 
(\ref{Fn}). Denote 
\begin{equation} \label{cR} {\mathcal R}=R\circ\Phi_5 
\end{equation} ${\mathcal R}$ is the germ of a function of 
five variables at $0$. We call it {\it the recovering 
invariant of} (\ref{affcs2}) {\it at the point $q_0$}. By 
construction, it is invariant of the weight $0$.  
 Note that 
the set of germs of systems of the type (\ref{affcs2}) 
having the reductions, which satisfies conditions 
(\ref{geng15}), (\ref{geng25}), and (\ref{geng35}), is 
generic. Using Theorem \ref{classn} in the case $n=5$ and 
the definition of the recovering invariant, we obtain the 
following classification of generic real analytic germs of 
systems of the type (\ref{affcs2}) in terms of the tuple of 
the primary invariants of their reductions and their 
recovering invariant: 
 
\begin{theor}
\label{class25} 
Given $4$ germs ${\mathcal I}_1$, ${\mathcal I}_2$, 
${\mathcal I}_3$, and ${\mathcal R}$ of real analytic 
functions of $5$ variables at $0$ such that ${\mathcal 
I}_3(0)\neq 0$, $10$ germs $\{\beta_{k}\}_{k=1}^5$ and 
$\{\phi_{k22}\}_{k=1}^5$ of real analytic functions of $4$ 
variables at $0$, $10$ germs $\{\psi_{k}\}_{k=1}^5$ and 
$\{\phi_{k23}\}_{k=1}^5$ of real analytic functions of $3$ 
variables at $0$, and $5$ germs $\{\phi_{k24}\}_{k=1}^5$ of 
real analytic functions of $2$ variables at $0$ there 
exists a unique, up to state-feedback real analytic 
transformation of the type (\ref{feedback}), 
five-dimensional real analytic control-affine system with 
two input such that first its reduction satisfies 
genericity assumptions (\ref{geng15}),(\ref{geng25}), 
(\ref{geng35}), secondly $\Bigl(\{{\mathcal 
I}_j\}_{j=1}^{3},\{\beta_{k}\}_{k=1}^5,\{\psi_{k}\}_{k=1}^5,\, 
\{\phi_{k2l}:1\leq k\leq 5, 2\leq l \leq 4\}\Bigr)$ is its 
tuple of the primary invariants, and finally ${\mathcal R}$ 
is its recovering invariant at the given point $q_0$. In 
other words, the tuples of the primary invariants give the 
regular $(3,2)$-parameterization of the considered 
classification problem in the real analytic category with 
the following $4\times 4$ parameterization matrix $P$: 
\begin{equation}
\label{parmat25} P=\begin{pmatrix} 0&0&0&5\\ 
0&0&0&10\\0&0&0&10\\1&0&0&3\end{pmatrix}. \end{equation} 
The Poincare series $M(t)$ of the considered classification 
problem satisfies 
\begin{equation}
\label{poin25} M(t)=\cfrac{1}{(1-t)^6}+t^3 
\left(\cfrac{5}{(1-t)^3}+\cfrac{10}{(1-t)^4}+\cfrac{10}{(1-t)^5}+
\cfrac{3}{(1-t)^6}\right). 
\end{equation}
\end{theor} 

By the same arguments, as in the case of $(1,n)$-affine 
control system with $n\geq 5$, treated in section 3, one 
can show that (3,2) is the characteristic pair of the 
considered classification problem (just use formulas 
(\ref{Nicein2}), (\ref{Nicein3}), and Corollary 
\ref{usecor}). The characteristic matrix ${\mathcal C}$ of 
the problem is equal to ${\rm Norm}(P)$, which can be found 
easily by the series of elementary transformations. Namely, 
\begin{equation}
\label{C52} {\mathcal C}={\rm Norm}(P)=\begin{pmatrix} 
1&3&6&5\\ 0&0&0&16\\0&0&0&13\\0&0&0&4\end{pmatrix}. 
\end{equation}
\begin{concl}
\label{concl25} The characteristic parameterization of 
$(2,5)$ control-affine systems, up to state-feedback 
transformations, consists of $4$ functional invariants of 
$5$ variables and the weight $3$, $13$ functional 
invariants of $4$ variables and the weight $3$, $16$ 
functional invariants of $3$ variables and the weight $3$, 
$1$ functional invariants of $2$ variables and the weight 
$0$, $3$ functional invariants of $2$ variables and the 
weight $1$, $6$ functional invariants of $2$ variables and 
the weight $2$, and $5$ functional invariants of $2$ 
variables and the weight $3$. 
\end{concl} 

In order to obtain a characteristic parameterization from 
the parameterization by the tuple of the primary invariants 
and the recovering invariant one can implement some series 
of rearrangement of the primary invariants according to the 
series of elementary transformations from the matrix $P$ to 
${\rm Norm}(P)$, as was described in the previous section 
(see, for example, formula (\ref{FG}) and the paragraph 
after it).

\end{document}